\documentclass[12pt]{amsart}

\usepackage{mathtools}
\mathtoolsset{showonlyrefs}
\usepackage{amssymb}

\usepackage{dirtytalk}

\usepackage{enumitem}
\usepackage[top=1.5in, bottom=1.5in, left=1.25in, right=1.25in]{geometry}
\usepackage[bottom]{footmisc}
\usepackage{fancyhdr}

\usepackage[backend=biber,style=numeric,sorting=nyt,maxbibnames=99,url=false,giveninits=true]{biblatex}
\AtEveryBibitem{
  \clearfield{pages}   
}
\usepackage{xcolor}
\usepackage{hyperref}
\usepackage{cleveref}

\hypersetup{
    colorlinks=true,
    linkcolor=blue,
    citecolor=magenta,
    urlcolor=cyan
}

\allowdisplaybreaks

\AtBeginDocument{\setlength{\footskip}{30pt}}
\newcommand{\abs}[1]{\left\lvert#1\right\rvert}
\newcommand{\comp}[1]{{#1}^{\mathsf{c}}}
\newcommand{\eps}{\varepsilon}
\newcommand{\norm}[1]{\left\lVert#1\right\rVert}
\newcommand{\one}{\mathbf{1}}
\newcommand{\set}[1]{\left\{#1\right\}}

\DeclareMathOperator{\dd}{\mathrm{d}}
\DeclareMathOperator{\dist}{\mathrm{dist}}
\DeclareMathOperator{\supp}{\mathrm{supp}}

\newcommand{\roverline}[1]{\mathpalette\doroverline{#1}}
\newcommand{\doroverline}[2]{\overline{#1#2}}

\newcommand{\C}{\mathbb{C}}
\newcommand{\D}{\mathbb{D}}
\newcommand{\E}{\mathbb{E}}
\newcommand{\F}{\mathbb{F}}
\newcommand{\N}{\mathbb{N}}
\renewcommand{\P}{\mathbb{P}}
\newcommand{\Q}{\mathbb{Q}}
\newcommand{\R}{\mathbb{R}}
\renewcommand{\S}{\mathbb{S}}
\newcommand{\T}{\mathbb{T}}
\newcommand{\Z}{\mathbb{Z}}

\newcommand{\BB}{\mathcal{B}}
\newcommand{\CC}{\mathcal{C}}
\newcommand{\EE}{\mathcal{E}}
\newcommand{\FF}{\mathcal{F}}
\newcommand{\II}{\mathcal{I}}
\newcommand{\RR}{\mathcal{R}}
\newcommand{\WW}{\mathcal{W}}

\newcommand{\0}{\mathbf{0}}

\newtheorem{theorem}{Theorem}[section]
\newtheorem{proposition}[theorem]{Proposition}
\newtheorem{cor}[theorem]{Corollary}
\newtheorem{lemma}[theorem]{Lemma}

\theoremstyle{remark}
\newtheorem{remark}[theorem]{Remark}
\newtheorem{ex}[theorem]{Example}

\AddToHook{env/proposition/begin}{\crefalias{theorem}{proposition}}
\AddToHook{env/lemma/begin}{\crefalias{theorem}{lemma}}
\AddToHook{env/remark/begin}{\crefalias{theorem}{remark}}

\crefname{cor}{corollary}{corollaries}
\Crefname{cor}{Corollary}{Corollaries}

\expandafter\patchcmd\csname\string\proof\endcsname
  {\itshape}
  {\color{purple}\itshape}
  {}{}

\newcommand{\authoraddress}[3]{%
  \address{%
    {\normalfont\bfseries #1}\endgraf
    {\normalfont\scshape #2}\endgraf
    {\normalfont\ttfamily #3}%
  }%
}

\begin{document}
\title{Gowers Uniformity and Pointwise (Non-)Convergence of Weighted Ergodic Averages}
\makeatletter
\let\inserttitle\@title
\makeatother

\author{F\'elix Brokering Pinilla}
\authoraddress
  {F\'elix Brokering Pinilla}
  {School of Mathematics, University of Bristol}
  {felix.brokeringpinilla@bristol.ac.uk}

\author{Ben Krause}
\authoraddress
    {Ben Krause}
    {School of Mathematics, University of Bristol}
    {ben.krause@bristol.ac.uk}


\fancypagestyle{myfancy}{%
  \fancyhf{}%
  \fancyhead[CE]{\scshape Uniformity of Weighted Ergodic Averages}%
  \fancyhead[CO]{F.\ Brokering Pinilla and B.\ Krause}%
  \fancyhead[LE,RO]{\thepage}%
  \renewcommand{\headrulewidth}{0pt}%
  \renewcommand{\footrulewidth}{0pt}%
}

\begin{abstract}

We study weighted ergodic averages of the form \[\frac{1}{N}\sum_{n=1}^{N}w(n)f(T^nx),\] where $(X,\mu,T)$ is a measure preserving system, $f\in L^{\infty}(\mu)$ and $w\colon \Z\to \C$ is a bounded weight. If for some $s\geq 2$, $\eps>0$ and a finite constant $C_{\eps,s}>0$ \[\norm{w}_{U^s[N]} \leq C_{\eps,s}(\log N)^{-1/2^s-\eps}\] for all $N\geq 2$, then the averages above converge almost surely; we prove that this result is sharp by constructing a $1$-bounded weight whose Gowers norms are all simultaneously bounded by $(\log N)^{-1/2^s}$ up to an arbitrarily slowly growing factor, but whose weighted ergodic averages oscillate between $-1$ and $1$ almost everywhere in a suitable system. In fact, we show that the property of failure of pointwise convergence is highly \say{generic} among bounded weights. The core idea behind the construction is to use multiplicative characters to build \textit{character words} which preserve earlier patterns on long blocks while gaining strong cancellation at every order. 

\end{abstract}

\maketitle
 
\pagestyle{myfancy}
\setcounter{tocdepth}{1}
\tableofcontents

\section{Introduction}

\subsection{History}

Weighted ergodic averages appear as one of the earliest variants of the landmark Ergodic Theorems of J. von Neumann \cite{VonNeumman} and G.D. Birkhoff \cite{Birkhoff}, dating back at least to the work of N. Wiener and A. Wintner \cite{WienerWintner}, who showed that for every measure preserving system $(X,\mu,T)$ and every $f\in L^1(\mu)$ there exists a full measure set on which the averages \[\frac{1}{N}\sum_{n=1}^{N}e(n\beta)f(T^nx), \qquad e(t) := e^{2\pi i t}\] converge for every $\beta \in \T$.

In full generality, the central issue is the norm and pointwise convergence of the averages 
\begin{equation}
\label{e:WeightedAverages}
    \frac{1}{N}\sum_{n=1}^{N}w(n)f(T^nx),
\end{equation}
where $(X,\mu,T)$ is a measure preserving system, $f$ is a measurable function on $X$ of suitable integrability and $w\colon \Z\to \C$ is a weight. We say that $w$ is \textit{norm good} if the averages \eqref{e:WeightedAverages} converge in $L^{2}$-norm for every measure preserving system $(X,\mu,T)$ and every $f\in L^{2}(\mu)$, and say that it is \textit{pointwise good for $L^{p}$ functions} for $p\geq 1$ if they converge pointwise almost everywhere whenever $f\in L^{p}(\mu)$. We are mostly interested in the case when $w$ is bounded. Here, the trivial estimate \[\abs{\frac{1}{N}\sum_{n=1}^{N}w(n)f(T^nx)} \leq \norm{w}_{\ell^{\infty}(\Z)}\frac{1}{N}\sum_{n=1}^{N}\abs{f}(T^nx)\] along with the boundedness of the ergodic maximal function allows us to pass pointwise convergence from functions in $L^{\infty}$ (or any dense subclass) to $L^{p}$ for every $p\geq 1$, so we may simply say that $w$ is \textit{pointwise good} and \textit{not pointwise good}\footnote{Although \textit{pointwise bad} would perhaps be a more convenient negation to pointwise good, it is already broadly used to mean that the failure of pointwise convergence happens on every aperiodic system. We stick with this slightly peskier naming convention hoping not to confuse the reader.} otherwise. By the Dominated Convergence Theorem, any bounded pointwise good weight is immediately norm good.

Early examples of pointwise good weights were found by C. Ryll-Nardzewski \cite{RyllNardzewski}, who showed that bounded weights that can be approximated by trigonometric polynomials in the Besicovitch distance, recalled in \eqref{e:DefiBesicovitch} below, are pointwise good. An important part of the foundations of the theory was then laid out by A. Bellow and V. Losert in their systematic study of weighted and subsequential ergodic theorems \cite{BellowLosert}, where they notably extended Ryll-Nardzewski's result in the spirit of Wiener and Wintner. Many additional examples also arise as generalizations of Wiener-Wintner, with an early extension to polynomial phases due to E. Lesigne \cite{LesignePolynomial} and a more recent far-reaching generalization to nilsequences due to B. Host and B. Kra \cite{HostKra}.

A further classical source of pointwise good weights is J. Bourgain's celebrated Return Times Theorem \cite{BourgainReturnTimes} (see also \cite{BFKO}), that states that for every measure preserving system $(Y,\nu,S)$ and every $g\in L^{\infty}(\nu)$, the sequence $w(n)=g(S^ny)$ is pointwise good for $\nu$-almost every $y\in Y$. This was then extended in one direction by D.J. Rudolph \cite{Rudolph}, who established pointwise goodness for multiterm return times of the form \[w(n)=g_1(S_1^{n}y_1)g_{2}(S_2^{n}y_2)\ldots g_k(S_k^{n}y_k),\] where each $g_i$ is a bounded function on some measure preserving system $(Y_i,\nu_i,S_i)$, $1\leq i \leq k$, and the full measure set for $y_{i}\in Y_{i}$ can be chosen before the functions $g_{i+1},\ldots,g_k$ and the points $y_{i+1},\ldots,y_k$ have been fixed. In a separate direction, P. Zorin-Kranich \cite{ZorinKranichDouble} showed that double return times \[w(n)=g_1(S^{an}y)g_{2}(S^{bn}y)\] are pointwise good for every measure preserving system $(Y,\nu,S)$,  $g_1,g_2\in L^{\infty}(\nu)$, $a,b\in \Z$, and almost every $y\in Y$.

In the arithmetic setting, the pointwise goodness for functions in $L^{p}$, $p>1$ of the von Mangoldt function follows by summation by parts from the pointwise ergodic theorem along the primes, which was established for $p> (1+\sqrt{3})/2$ by Bourgain \cite{BourgainPrimes} and extended to $p>1$ by M. Wierdl \cite{WierdlPrimes}. For the M\"obius and Liouville functions we have pointwise goodness in the full range, as shown by E.H. El Abdalaoui, J. Kułaga-Przymus,
M. Lemańczyk and T. de la Rue \cite{EKLR}, and in fact the same is true in the broader class of $1$-bounded multiplicative functions that have a mean on every arithmetic progression, as observed by N. Frantzikinakis and Host \cite{FrantzikinakisHost2017}.

For more examples and a complete account on weighted ergodic theorems and return times theorems we refer to the dedicated chapter in \cite{EisnerFarkas}.\\

Crucially, each of the results above makes use of the concrete structure of the weights it treats, however broad that class may be. A key role in giving general criteria for norm or pointwise goodness is played by the Spectral Theorem, which states that for every measure preserving system $(X,\mu,T)$, and every $f\in L^{2}(\mu)$, there exists a finite Borel measure $\sigma_{f}$ on $\T$ so that \[\norm{\frac{1}{N}\sum_{n=1}^{N}w(n)T^nf}_{L^2(\mu)} = \norm{\frac{1}{N}\sum_{n=1}^{N}w(n)e(n\beta)}_{L^{2}(\sigma_{f})}\] for every $w\colon \Z\to \C$ and every $N\geq 1$. This readily characterizes norm goodness, as one can see that a bounded weight $w$ is norm good if and only if the weighted exponential averages
\begin{equation}
\label{e:ExponentialAverage}
    \frac{1}{N}\sum_{n=1}^{N}w(n)e(n\beta)
\end{equation}
converge for every $\beta\in \T$.

Pointwise goodness inevitably requires a bit more. After possibly centering $w$, the quantity to look at is 
\begin{equation}
\label{e:LargestFourier}
    \sup_{\beta\in \T}\abs{\frac{1}{N}\sum_{n=1}^{N}w(n)e(n\beta)},
\end{equation}
for $N\geq 1$, which captures the uniform rate of convergence of the averages \eqref{e:ExponentialAverage} to their limit; slightly abusing notation, we call \eqref{e:LargestFourier} the \textit{largest (normalized) Fourier coefficient} of $w$. If a weight exhibits sufficiently fast decay in this quantity, which we will simply call \textit{Fourier decay}, then it will be pointwise good. The following result was first formulated for the M\"obius function and $\eta>1$ by El Abdalaoui, Kułaga-Przymus, Lemańczyk and de la Rue \cite{EKLR} and later refined to $\eta>1/2$ and any bounded weight $w$ by A.-H. Fan \cite{Fan}, although it already appears implicitly in work of Bourgain \cite{B0}. We refer the reader to \Cref{ss:Not} to recall the $\lesssim$ notation.

\begin{theorem}
\label{t:PositiveResultFourier}
    Let $w\colon \Z\to \C$ be a bounded weight. Suppose there exists some $\eta > 1/2$ such that for every $N\geq 2$ \[\sup_{\beta\in \T}\abs{\frac{1}{N}\sum_{n=1}^{N}w(n)e(n\beta)}\lesssim \left(\frac{1}{\log N}\right)^{\eta}.\] Then $w$ is pointwise good.
\end{theorem}

The core of this result is the elemental inequality \[\norm{\frac{1}{N}\sum_{n=1}^{N}w(n)T^nf}_{L^{2}(\mu)} \leq \sup_{\beta\in \T}\abs{\frac{1}{N}\sum_{n=1}^{N}w(n)e(n\beta)} \norm{f}_{L^2(\mu)},\] which follows from the Spectral Theorem. The threshold of $1/2$ for the exponent $\eta$ shows up naturally in the proof that rests on a Borel-Cantelli argument in $L^{2}$, which is where the spectral theory forces us to work.

\subsection{Main Results}

To improve on this exponent, we propose a parallel estimate where we trade off decay in higher $L^{p}$-norms of the averages for higher uniformity of $w$, as measured by Gowers norms, whose definition we recall in \Cref{sec:Prelims}. These norms measure the correlation of $w$ with sequences of higher complexity and have played a crucial role in the study of weighted ergodic averages, starting with work by Frantzikinakis, Host and Kra to control prime-weighted averages in norm \cite{FHK2007,FHK2013}, with further developments in the pointwise setting by J. Ter\"av\"ainen  \cite{Teravainen}, J. Fornal, A. Fragkos, M. Lacey, H. Mousavi, Y.-C. Sun and the second author \cite{FornalWW} and Fornal and the second author \cite{FornalKrause2026}. In our case, the relevant estimate for the averages \eqref{e:WeightedAverages} is the following.

\begin{lemma}
\label{l:L2SUS}
    Let $w\colon \Z\to \C$ and let $s\geq 2$. Then, for every measure preserving system $(X,\mu,T)$ and $f\in L^{2^s}(\mu)$, we have 
    \begin{equation}
    \label{e:Usestimate}
        \norm{\frac{1}{N}\sum_{n=1}^{N}w(n)T^nf}_{L^{2^s}(\mu)} \lesssim_{s} \norm{w}_{U^s[N]}\norm{f}_{L^{2^s}(\mu)}
    \end{equation}
    for every $N\geq 1$.
\end{lemma}

The following positive result then follows from \Cref{l:L2SUS} and a Borel-Cantelli argument, after a standard passage to lacunary times.
\begin{proposition}
\label{p:PositiveResult}
    Let $w\colon \Z\to\C$ be a bounded weight and suppose that there exist some $s\geq 2$ and $\eta > \frac{1}{2^s}$ so that \[\norm{w}_{U^s[N]} \lesssim_{s} \left(\frac{1}{\log N}\right)^{\eta},\] for every $N\geq 2$. Then $w$ is pointwise good.
\end{proposition}

While the proof of \Cref{l:L2SUS} uses fairly elementary tools ---Calder\'{o}n's Transference Principle, the Hausdorff-Young inequality, and induction--- our main result shows that it is essentially sharp.

\begin{theorem}
\label{t:MainResult}
    For any non-decreasing function $\psi\colon \R_{+}\to \R_{+}$ with $\psi(x) \to \infty$ as $x\to \infty$, there exists a weight $w\colon \Z \to \set{-1,0,1}$ such that for every $N\geq 2$ \[\sup_{\beta\in \T}\abs{\frac{1}{N}\sum_{n=1}^{N}w(n)e(n\beta)} \lesssim \left(\frac{\psi(\log N)}{\log N}\right)^{1/2},\] and \[\norm{w}_{U^s[N]} \lesssim_{s} \left(\frac{\psi(\log N)}{\log N}\right)^{1/2^s},\] simultaneously for every $s\geq 2$, and a measure preserving system $(X,\mu,T)$ along with a function $f\colon X\to \set{-1,1}$ so that  \[\liminf_{N\to\infty}\frac{1}{N}\sum_{n=1}^{N}w(n)f(T^nx) = -1, \qquad \limsup_{N\to\infty}\frac{1}{N}\sum_{n=1}^{N}w(n)f(T^nx) = 1\] pointwise almost everywhere.
\end{theorem}

Since $\psi$ can be taken to grow arbitrarily slowly, we can approach the thresholds with arbitrary precision. For instance, choosing $\psi(n)= \log n$ we can have $w$ satisfy \[\sup_{\beta\in \T}\abs{\frac{1}{N}\sum_{n=1}^{N}w(n)e(n\beta)} \lesssim_{\eps} \left(\frac{1}{\log N}\right)^{1/2-\eps}, \quad \norm{w}_{U^s[N]} \lesssim_{s,\eps} \left(\frac{1}{\log N}\right)^{1/2^s-\eps},\] simultaneously for every $s\geq 2$ and $\eps>0$, while still having the averages \eqref{e:WeightedAverages} exhibit the full oscillation permitted by its bounds. The exact endpoint where $\psi$ is absent is left open by our methods.\\

An important takeaway from our construction is that it shows that norm and pointwise convergence are quantitatively very different phenomena. To make this more precise, we need another definition. We call a weight $w\colon \Z\to \C$ \textit{ergodic} if for every ergodic measure preserving system $(X,\mu,T)$ and $f\in L^{2}(\mu)$ the $L^{2}$-limit of the weighted ergodic averages \eqref{e:WeightedAverages} is the integral of $f$. We will say that $w\colon \Z\to\C$ is \textit{strongly ergodic} if for every $k,\ell\geq 1$ the weighted multiple ergodic averages
\begin{equation}
\label{e:WeightedMultipleAverages}
    \frac{1}{N}\sum_{n=1}^{N}w(n)f_1\left(\prod_{i=1}^{\ell}T_{i}^{p_{i,1}(n)}x\right)\cdot \ldots \cdot f_k\left(\prod_{i=1}^{\ell}T_{i}^{p_{i,k}(n)}x\right),
\end{equation}
converge in $L^2$-norm to the same $L^{2}$-limit as the unweighted averages (which is known to exist by M. N. Walsh's Theorem \cite{Walsh}), where 
\begin{itemize}
    \item $T_1,\ldots,T_{\ell}$ are measure preserving transformations on a probability space $(X,\mu)$ generating a nilpotent group;
    \item $f_1,\dots,f_k \in L^\infty(\mu)$;
    \item $p_{i,j}\colon \Z\to \Z$ are integer polynomials for  $1 \leq i \leq \ell, \ 1 \leq j \leq k$.
\end{itemize}
We say that $w$ is \textit{strongly null} if the $L^{2}$-limit of the averages \eqref{e:WeightedMultipleAverages} is always $0$. Any strongly ergodic weight is ergodic, and both ergodic and strongly null weights are norm good. Moreover, by A. Leibman's nilpotent PET induction scheme \cite{Leibman} (see also \cite{FHK2013, FrantzikinakisHost2017}), we can bound the $L^2$-norm of \eqref{e:WeightedMultipleAverages} by a sufficiently large Gowers norm of $w$, so any weight with $\norm{w}_{U^s[N]}\to 0$ as $N\to \infty$ for every $s\geq 2$ is strongly null.

Choosing $\psi$ so that the Gowers norms of the weight given by \Cref{t:MainResult} decay to $0$ yields a strongly null weight that fails to be pointwise good. Adding $1$ to the weight immediately yields the following sharp qualitative divide between norm and pointwise convergence.

\begin{cor}
\label{c:StronglyErgodicNotPointwiseGood}
    There exists a weight $w\colon \Z\to \set{0,1,2}$ that is strongly ergodic and not pointwise good.
\end{cor}

Most examples of weights which are norm good but not pointwise good that we are aware of are either unbounded or fail to have decaying largest Fourier coefficients (and Gowers norms) after centering. The only two exceptions we know of arise as implicit consequences of an example due to J.-P. Thouvenot and B. Weiss reported in \cite{LinOlsenTempelman} and described in unpublished notes \cite{ThouvenotWeiss} that Weiss kindly shared with us, and a construction of M. Hochman that was communicated to us by F. Richter and can be found in the appendix of \cite{Tal}. In both cases, they obtain the existence of a weakly mixing system $(X,\mu,T)$ and a generic point $x\in X$ such that the weight $w(n)=f(T^nx)$ fails to be pointwise good for some continuous function $f\colon X\to \C$, which after centering can be taken to have mean zero. Since $(X,\mu,T)$ is weakly mixing, all Host-Kra seminorms of $f$ vanish, and since $x$ is generic, this implies that $\norm{w}_{U^s[N]} \to 0$ as $N\to \infty$ for every $s\geq 2$. This means $w$ is in fact strongly null, so both examples readily give the same initial qualitative conclusion as our result.

Our construction improves on them by providing explicit examples with concrete decay rates. In contrast, Hochman proves the existence of an example by a Baire category argument, so no quantitative decay bounds are available, and while the explicit construction of Thouvenot and Weiss can be made quantitative, the rates one obtains do not depend directly on $N$.

\subsection{Overview of the Construction}

Let us begin by describing exactly what \Cref{t:MainResult} must accomplish.

On the one side, we need strong cancellation: the largest Fourier coefficient and the Gowers norms of $w$ all have to vanish as $N\to \infty$. This should be understood as a randomness property, since it requires the weights to lack correlations with structured sequences. As observed above, it forces the averages \eqref{e:WeightedAverages} to converge to $0$ in $L^2$-norm.

On the other hand, $w$ must be sufficiently structured to persistently correlate with an orbit: since the $L^2$-limit of the weighted averages is zero, pointwise failure requires us to find a measure preserving system $(X,\mu,T)$ and a bounded function $f$ so that \[\limsup_{N\to \infty} \abs{\frac{1}{N} \sum_{n=1}^{N}w(n)f(T^nx)} > 0\] for $x$ in a set of positive measure. This can be regarded as a self-similarity property; since sets of positive measure are visited infinitely often, the above abstracts to finding a bounded sequence, $a\colon \Z\to \C$, so that 
\begin{equation}
\label{e:SequenceCorrelations}
    \limsup_{N\to \infty} \abs{\frac{1}{N} \sum_{n=1}^{N}w(n)a(n+m)} > 0
\end{equation}
for many shifts $m\geq 1$.\\

These two features evidently pull in opposite directions. In particular, they rule out the usual examples of traveling blocks which are constant on increasingly long intervals; they are self-similar enough to guarantee pointwise non-convergence, but give very little room for cancellation. Our starting point to overcome this difficulty is an example proposed by Bourgain in the appendix of his unpublished notes on Return Times \cite{BourgainUnpublished}. Here, he proposes replacing the constants in the traveling blocks by wave packets \[\sum_{c\in I}\xi_{c}e(n\alpha_c)\] with amplitudes $\xi_c\in \C$ and frequencies $\alpha_c\in \T$ varying between intervals. 

The construction in these unpublished notes is, however, incomplete, and one can check that it fails to have decaying largest Fourier coefficients, not to mention Gowers norms. It is in fact not even norm good. Nevertheless, developing his approach led us to construct a first example of a weight with concrete quantitative decay rates that fails to be pointwise good. 

To achieve a sharp result, we passed from the Fourier-analytic framework of Bourgain to an arithmetic one, where destructive interference can be most sharply captured. In particular, we pass to a finite-field setting and replace exponential phases by multiplicative characters, transferring estimates to the integers while at the same time retaining self-similarity across scales by introducing \textit{character words}. This type of argument goes back at least to work of M. Christ on ergodic averages over sparse sequences \cite{Christ}, where a similar problem of transferring bounds from finite fields to $\Z$ is solved by first passing to a suitable larger cyclic group. It also takes inspiration in the work of M. Keane on generalized Morse sequences \cite{Keane}, where a notion of product of words is introduced that allows the construction of sequences with prescribed self-similarity.

Roughly speaking, a character word is a word on the symbols $\set{-1,1}$ constructed by choosing some word $v\in \set{-1,1}^{d}$ and an odd prime $p$, and concatenating $p$ copies of $v$ while changing the signs of the $i$-th copy according to the value of quadratic character of $\F_p$ on $i$. This results in a word $V\in \set{-1,1}^{dp}$ that locally looks like $v$, while having the cancellation properties of the quadratic character of $\F_p$. Concretely, we show that if we regard $V$ as a function on $\Z_{\ell}$, where $\ell=dp$, we have the optimal estimates \[\sup_{\xi\in \Z_{\ell}} \abs{\widehat{V}(\xi)} \lesssim \frac{1}{p^{1/2}}, \quad \norm{V}_{U^s(\Z_{\ell})} \lesssim_{s} \frac{1}{p^{1/2^s}},\ \forall s\geq 2.\] 

To pass all the way to the integers, we construct a sequence of nested character words by using each word to define the next one. We then define a weight in the same spirit as above by replacing wave packets by character words on each interval. Choosing the primes that define the character words and the lengths of the intervals appropriately allows us to transfer the desired cancellation to each scale. To pick up correlations, we observe that any translate of a character word will locally look like a signed translate of the preceding ones, so the correlation of the weight with a shift of itself will be comparable to the correlation of two signed translates of a common character word; we make this quantity large by simply letting the weight run through each possible signed translate of every character word instead of having each character word appear just once. The quantitative gain in comparison with the construction proposed by Bourgain comes from the lower-dimensionality of the class of character words versus suitable wave packets.

\subsection{Further Results}

As a consequence of our construction, we show that the property of not being pointwise good is generic for bounded weights, even for those with prescribed Fourier and Gowers decay arbitrarily close to the thresholds of \Cref{t:PositiveResultFourier} and \Cref{p:PositiveResult}. Concretely, we prove that after endowing the space of bounded weights with the Besicovitch distance,
\begin{equation}
\label{e:DefiBesicovitch}
    d_B(u,v)\coloneqq \limsup_{N\to\infty}\frac{1}{N}\sum_{n=1}^{N}\abs{u(n)-v(n)},
\end{equation}
and identifying weights at distance zero, the set of weights that are not pointwise good is open and dense, even when restricting to the subspace of weights satisfying the decay rates of \Cref{t:MainResult} for a fixed $\psi$; these are the contents of \Cref{t:WeightGeneric} and \Cref{c:PrescribedDecayGeneric}, respectively.

The above means that the property of not being pointwise good is very easy to achieve ---arbitrarily small perturbations in mean suffice to destroy pointwise goodness--- and yet rather hard to get rid of ---even changes in positive proportion may not fix it. In particular, it stands in striking contrast to the many examples discussed above, and most specially with Bourgain’s Return Times Theorem, which alone supplies a huge class of pointwise good weights. Our result implies that they all nevertheless live inside of a closed nowhere dense set for $d_{B}$ and are arbitrarily close to not being pointwise good.\\

In a related direction, there is a close relationship between weighted and sub-sequential ergodic averages: for instance, summation by parts shows that for any bounded $a\colon \Z\to \C$,
\begin{align}
    \frac{1}{N} \sum_{n=1}^{N} a(n^2) \to L \iff \frac{1}{N} \sum_{n=1}^{N} w_{\square}(n) \cdot a(n) \to L, \quad w_{\square}(n) := 2\sqrt{n} \cdot \one_{\set{k^2 \colon k \geq 1}}(n).
\end{align}
Accordingly, our analysis also allows us to address the subsequential setting, specifically in the case of \textit{Gowers uniform} sets, which are particular subsets of the integers whose asymptotic statistics with respect to combinatorial patterns of finite complexity are determined exactly by their density; for a precise definition, examples and a short discussion of their properties we refer to \Cref{sec:GowersUniform}. Crucially, these sets enjoy strong recurrence and norm convergence properties, as shown by Frantzikinakis and Host \cite{FrantzikinakisHost2017}. Nevertheless, combining our construction with Bourgain's probabilistic balayage technique \cite[Appendix~3]{BourgainUnpublished}, we are able to obtain arbitrarily
dense Gowers uniform sets along which the
pointwise ergodic theorem fails. 

To make this precise, we recycle nomenclature from the weighted setting, and call a subset $E\subseteq \N$ \textit{norm good} or \textit{pointwise good for $L^p$ functions} for $p\geq 1$, if for every measure preserving system $(X,\mu,T)$ the averages
\begin{equation}
\label{e:SubseqAverages}
    \frac{1}{\abs{E\cap [N]}}\sum_{n\in E\cap [N]}f(T^nx)
\end{equation}
converge in $L^{2}(\mu)$ for $f\in L^{2}(\mu)$ or pointwise almost everywhere for $f\in L^{p}(\mu)$, respectively. If $E$ admits a positive density, we can rewrite \eqref{e:SubseqAverages} in terms of the ergodic average weighted by $\one_{E}$ and all pointwise goodness notions coincide. What we show is the following.

\begin{theorem}
\label{t:GowersUniformBad}
    For every $0<\eps\leq 1/2$, there exists a set $E\subseteq \N$ with $d(E)= 1-\eps$ that is Gowers uniform and not pointwise good.
\end{theorem}
Density 1 itself is of course different, since every such set is immediately pointwise good. As in our main theorem, the conclusion can also be phrased in a quantitative way and is again generic for the natural distance in the appropriate space of subsets of $\N$.\\

As a final application, our mechanism can be adapted to produce strongly ergodic weights, with the same quantitative bounds as in \Cref{t:MainResult} that nevertheless exhibit the \textit{strong sweeping out property}, which is the most drastic form of pointwise non-convergence. Since this phenomenon concerns unbounded weights and this paper is mainly focused on bounded weights, we carry this out in \Cref{appendix:Sweeping}.

\subsection{Organization of the Paper} 
The structure of the paper is as follows:
\begin{itemize}
    \item We begin by consolidating the required machinery and proving \Cref{p:PositiveResult} in \Cref{sec:Prelims};
    \item \Cref{sec:CharacterWord} then develops character words and establishes the necessary uniformity estimates;
    \item The construction is carried out in \Cref{sec:Construction}, where we prove \Cref{t:MainResult};
    \item The genericity results are handled in \Cref{sec:Genericity};
    \item We address Gowers uniform sets in \Cref{sec:GowersUniform}, where we prove \Cref{t:GowersUniformBad};
    \item We discuss the extension to weights with the strong sweeping-out property in \Cref{appendix:Sweeping}.
\end{itemize}
\subsection*{Acknowledgments}

We would like to express our gratitude to Yu-Chen Sun for the insightful suggestion of looking at multiplicative functions for a better source of cancellation, and for patiently explaining the basics. The first author also thanks Van Cyr and Bryna Kra for valuable questions and comments during a seminar talk on this work.

\subsection*{Use of AI}

The authors acknowledge the use of LLMs for assistance in suggesting relevant literature and reviewing preliminary drafts throughout the development of this work. In particular, an audit of an earlier draft suggested looking at generalized Morse sequences as a general framework for constructing self-similar sequences, which led to the development of character words. The contents of this manuscript are the work of the authors, who take full responsibility for all mathematical statements and proofs.

\section{Preliminaries}
\label{sec:Prelims}

\subsection{Notation and Conventions}\label{ss:Not}

A \textit{measure preserving system} is a quadruple $(X,\FF,\mu,T)$ where $(X,\FF,\mu)$ is a probability space and $T\colon X\to X$ is a measure preserving transformation, this is, \[\mu(T^{-1}A) = \mu(A), \text{ for every }A\in \FF.\] We usually drop the $\sigma$-algebra and simply write $(X,\mu,T)$. Given a measurable function $f$ on $X$, we put $Tf \coloneqq f\circ T$.

Index ranges such as $n\geq 1$ or $k\geq 0$ and averaging cutoffs are always integers. We let $\Z_{+}=\set{0,1,2,3,\ldots}$ and $\N = \set{1,2,3,\ldots}$ and for $N\geq 1$ write $[N]=\set{1,\ldots, N}$. Given a set $E$, we let $\one_{E}$ denote its indicator function, and if $E\subseteq \Z$ we call \[d(E) \coloneqq \lim_{N\to \infty} \frac{\abs{E \cap [N]}}{N} = \lim_{N\to \infty} \frac{1}{N}\sum_{n=1}^{N}\one_{E}(n)\] the \textit{density} of $E$, whenever the limit exists; the limit superior always exists and is called the  \textit{upper density} of $E$ and denoted by $\overline{d}(E)$. We let $\overline{\D}$ denote the closed unit disk $\set{z\in \C \colon \abs{z}\leq 1}$, $\S$ denote the unit circle $\set{z\in \C \colon \abs{z}=1}$ and $\T$ denote the torus $\R/\Z$, which we identify with the interval $[0,1)$ and endow with the norm $\norm{\beta}_{\T} = \min_{n\in \Z}\abs{\beta - n}$. We use the standard abbreviation $e(t) = e^{2\pi i t}$ for $t\in \R$. 

We make use of modified Vinogradov notation and given two quantities $X,Y\geq 0$ write $X\lesssim Y$ or $Y\gtrsim X$ as a shorthand for the inequality $X\leq CY$ for some implicit constant $C>0$. If both $X\lesssim Y$ and $Y\lesssim X$ hold, we write $X\approx Y$. We indicate the dependence of this constant on a parameter $p$ by subscripts $\lesssim_{p}$, $\gtrsim_{p}, \approx_{p}$, and in particular, whenever $X$ and $Y$ are functions of more than one variable, the implicit constants in statements of the form \say{$X(p,q)\lesssim_{p}  Y(p,q)$ for every $q$} are understood to be uniform in $q$.

\subsection{Gowers Norms}

We employ Gowers norms in several settings. First, for a finitely supported function $a\colon \Z\to \C$ and $s\geq 1$ we define the \textit{unnormalized Gowers norms} via \[\norm{a}_{U^s(\Z)}^{2^s} = \sum_{n,h_1,\ldots,h_s\in \Z} \prod_{\omega \in \set{0,1}^{s}}\CC^{\abs{\omega}}a(n+\omega\cdot h),\] where $\CC$ denotes complex conjugation $\CC(z)=\overline{z}$, $\abs{\omega} = \omega_1 + \ldots + \omega_{s}$, $h=(h_1,\ldots,h_s)$ and $\omega\cdot h = \omega_1h_1+\ldots +\omega_{s}h_s$. These are norms for $s\geq 2$ and a seminorm for $s=1$. Moreover, they are translation, reflection and conjugation invariant for every $s\geq 1$ and additionally modulation invariant for $s\geq 2$. If we consider the \textit{discrete differencing operators} $\Delta_{h}$ defined for $h\in \Z$ by \[\Delta_{h}a(n) \coloneqq a(n)\overline{a(n+h)},\] we then have the following recursion for $s\geq 1$
\begin{equation}
\label{e:GowersZRecursion}
    \norm{a}_{U^{s+1}(\Z)}^{2^{s+1}} = \sum_{h\in \Z}\norm{\Delta_{h}a}_{U^s(\Z)}^{2^s},
\end{equation}
and, crucially, if we let \[\widehat{a}\colon \T\to \C, \qquad \widehat{a}(\beta) \coloneqq \sum_{n\in \Z}a(n)e(-n\beta)\] denote the Fourier Transform of $a$, we have the following identity 
\begin{equation}
\label{e:U2Fourier}
    \norm{a}_{U^2(\Z)} = \norm{\widehat{a}}_{L^4(\T)}.
\end{equation}

Starting from these norms we define the \textit{normalized Gowers norms} via \[\norm{a}_{U^s[N]}^{2^s} = \frac{1}{N^{s+1}}\norm{a\cdot \one_{[N]}}_{U^s(\Z)}^{2^s} = \frac{1}{N^{s+1}}\sum_{\substack{n,h_1,\cdots,h_s \in \Z \\ n+\omega\cdot h \in [N],\ \forall \omega \in \set{0,1}^{s}}}\prod_{\omega \in \set{0,1}^{s}}\CC^{\abs{\omega}}a(n+\omega\cdot h),\] for an integer $N\geq 1$. In the space of functions $a\colon [N]\to \C$ these are again modulation invariant norms for $s\geq 2$, whereas for $s=1$ we get the seminorm \[\norm{a}_{U^1[N]} = \abs{\frac{1}{N}\sum_{n=1}^{N}a(n)}.\] Normalized Gowers norms are increasing in the sense that $\norm{a}_{U^s[N]} \lesssim_{s} \norm{a}_{U^{s+1}[N]}$ and the recursion and Fourier identities above immediately yield the analogues \[\norm{a}_{U^{s+1}[N]}^{2^{s+1}} = \frac{1}{N}\sum_{h\in \Z}\norm{\Delta_{h}(a\cdot\one_{[N]})}_{U^s[N]}^{2^s}, \qquad \norm{a}_{U^2[N]} = \frac{1}{N^{3/4}}\norm{\widehat{a\cdot \one_{[N]}}}_{L^4(\T)}.\] Furthermore, the translation and conjugation invariance of the unnormalized norms imply the identity \[\norm{\Delta_{h}(a\cdot \one_{[N]})}_{U^s[N]} = \norm{\Delta_{-h}(a\cdot \one_{[N]})}_{U^s[N]}.\] 

More generally, we use the \textit{interval Gowers norms}, defined for a sequence $a\colon \Z\to \C$ and a non-empty finite interval $I\subseteq \Z$ as \[\norm{a}_{U^s[I]}^{2^s} = \frac{1}{\abs{I}^{s+1}}\norm{a\cdot \one_{I}}_{U^s(\Z)}^{2^s};\] note that $U^s[N]$-norms are precisely $U^s[I]$-norms for $I=[N]$. We only make use of the facts that translations of the input function correspond to translations of the interval \[\norm{a(\cdot+m)}_{U^s[I]} = \norm{a}_{U^s[I+m]},\ m\in \Z\] and that we have the trivial inequality \[\norm{a}_{U^s[I]} \leq \norm{a}_{\ell^{\infty}(I)}\] along with the identity \[\norm{a\cdot \one_{I}}_{U^s[N]} = \left(\frac{\abs{I}}{N}\right)^{\frac{s+1}{2^s}}\norm{a}_{U^s[I]}\] for intervals $I\subseteq [N]$.

Finally for functions defined on the cyclic groups $\Z_{\ell} = \Z/\ell \Z$, $\ell \geq 2$, we consider the \textit{cyclic Gowers norms} \[\norm{a}_{U^s(\Z_{\ell})}^{2^s} = \frac{1}{\ell^{s+1}}\sum_{n,h_1,\ldots,h_s\in \Z_{\ell}}\prod_{\omega \in \set{0,1}^{s}}\CC^{\abs{\omega}}a(n+\omega\cdot h).\] These are once again norms for $s\geq 2$ and only a seminorm for $s=1$, and remain invariant under translations, reflections and conjugations. For $s\geq 2$ they are also invariant under multiplications by additive characters of the corresponding group. The only additional property that we will need is the inequality \[\norm{a\cdot \one_{E}}_{U^s(\Z_{\ell})} \leq \left(\frac{\abs{E}}{\ell}\right)^{\frac{s+1}{2^s}}\norm{a}_{\ell^{\infty}(E)},\] for $E\subseteq \Z_{\ell}$. We observe that this definition makes sense also for $\ell$-periodic functions on $\Z$.

\subsection{Controlling Averages}

We now work our way towards \Cref{l:L2SUS}. The starting point is the following inequality in the integer model.

\begin{lemma}
\label{l:ell4-U2-control}
    Let $w,f\colon \Z\to \C$ be finitely supported. Then, \[\norm{\sum_{n\in \Z}w(n)f(x+n)}_{\ell^{4}(\Z)} \leq \norm{w}_{U^2(\Z)} \norm{f}_{\ell^{2}(\Z)}.\]
\end{lemma}

\begin{proof}
    By Hausdorff-Young we can bound \[\norm{\sum_{n\in \Z}w(n)f(x+n)}_{\ell^{4}(\Z)} = \norm{w(-\cdot) \ast f}_{\ell^4(\Z)} \leq \norm{\widehat{w(-\cdot)\ast f}}_{L^{4/3}(\T)} = \norm{\widehat{w(-\cdot)}\cdot \widehat{f}}_{L^{4/3}(\T)},\] and then by H\"older \[\norm{\sum_{n\in \Z}w(n)f(x+n)}_{\ell^{4}(\Z)} \leq \norm{\widehat{w(-\cdot)}}_{L^4(\T)} \norm{\widehat{f}}_{L^2(\T)}.\] Employing the Fourier identity for the $U^2(\Z)$-norm and Plancherel we then conclude \[\norm{\sum_{n\in \Z}w(n)f(x+n)}_{\ell^{4}(\Z)} \leq \norm{w(-\cdot)}_{U^2(\Z)}\norm{f}_{\ell^2(\Z)} = \norm{w}_{U^2(\Z)} \norm{f}_{\ell^2(\Z)},\] where we use that unnormalized Gowers norms are reflection invariant.
\end{proof}

Truncating $w \to w\one_{[N]}$ at each scale and employing Calder\'on's Transference Principle \cite{C1}, we immediately get the following ergodic consequence. We provide details for the unfamiliar reader.

\begin{lemma}
    For every measure preserving system $(X,\mu,T)$, $f\in L^{4}(\mu)$ and $w\colon \Z\to \C$ we have \[\norm{\frac{1}{N}\sum_{n=1}^{ N}w(n)T^nf}_{L^4(\mu)} \lesssim \norm{w}_{U^2[N]}\norm{f}_{L^4(\mu)}\] for every $N\geq 1$.
\end{lemma}

\begin{proof}
    Given $N\geq 1$ and $x\in X$ put \[F_{x}(n) = \begin{cases}
        f(T^nx),\ &\text{ if }1\leq n \leq 2N,\\
        0,\ &\text{ otherwise},
    \end{cases}\] so that for $1\leq k \leq N$ and $1\leq n \leq N$ we have $1\leq k+n \leq 2N$ and hence \[\sum_{n=1}^{N}w(n)F_{x}(k+n) = \sum_{n=1}^{N}w(n)f(T^{n}(T^kx)).\] By our estimate above and Cauchy-Schwarz we then have \[\sum_{k=1}^{N}\abs{\sum_{n=1}^{N}w(n)F_{x}(k+n)}^{4} \lesssim \norm{w}_{U^2(\Z)}^{4} N \sum_{k=1}^{ 2N}\abs{F_x(k)}^{4}.\] Integrating both sides with respect to $\mu$ we get \[\sum_{k=1}^{N}\int \abs{\sum_{n=1}^{N}w(n)f(T^n(T^kx))}^{4} \dd\mu \lesssim \norm{w}_{U^2(\Z)}^{4} N \sum_{k=1}^{2N}\int \abs{f(T^kx)}^{4}\dd\mu,\] and by invariance of the measure \[N \norm{\sum_{n=1}^{N}w(n)T^nf}_{L^4(\mu)}^{4} \lesssim \norm{w}_{U^2(\Z)}^{4}\cdot N^2 \norm{f}_{L^4(\mu)}^{4},\] so dividing by $N^{5}$ yields \[\norm{\frac{1}{N}\sum_{n=1}^{N}w(n)T^nf}_{L^4(\mu)}^{4} \lesssim \norm{w}_{U^2[N]}^{4}\norm{f}_{L^4(\mu)}^{4}\] and we are done.
\end{proof}

\begin{remark}
    Perhaps unsurprisingly, the estimate above is actually yet another manifestation of the Spectral Theorem. Indeed, considering the trivial $L^{2}(\T)\times L^{\infty}(\mu)\to L^{\infty}(\mu)$ bound given by Cauchy-Schwarz and Plancherel \[\norm{\frac{1}{N}\sum_{n=1}^{N}w(n)T^nf}_{L^{\infty}(\mu)} \leq N^{-1/2}\norm{\widehat{w\one_{[N]}}}_{L^2(\T)}\norm{f}_{L^{\infty}(\mu)}\] and the $L^{\infty}(\T)\times L^{2}(\mu)\to L^{2}(\mu)$ bound given by the Spectral Theorem \[\norm{\frac{1}{N}\sum_{n=1}^{N}w(n)T^nf}_{L^2(\mu)} \leq N^{-1} \norm{\widehat{w\one_{[N]}}}_{L^{\infty}(\T)}\norm{f}_{L^2(\mu)},\] one can obtain the estimate above by bilinear interpolation at the midpoint $L^{4}(\T) \times L^{4}(\mu) \to L^{4}(\mu)$ and the Fourier identity for $U^2[N].$
\end{remark}

Induction now yields \Cref{l:L2SUS}:

\begin{proof}[Proof of \Cref{l:L2SUS}]
    Assume the result holds for some $s\geq 2$ and any measure preserving system, function in $L^{2^s}$, and weight. Fixing $N\geq 1$, abbreviate $w=w\one_{[N]}$ and for $f\in L^{2^{s+1}(\mu)}$ pointwise expand \[\abs{\frac{1}{N}\sum_{n=1}^{N}w(n)T^nf}^2 = \frac{1}{N^2}\sum_{n=1}^{N}\abs{w(n)}^{2}\abs{T^nf}^2 + 2\text{Re}\left(\frac{1}{N^2}\sum_{h=1}^{N-1}\sum_{n=1}^{N-h}\Delta_{h}w(n)T^n(f\cdot \overline{T^hf})\right).\]
    Writing $\abs{w(n)}^2 = \Delta_{0}w(n)$ and extending the innermost sum in the second term to all of $[N]$ by recalling that $w=w\one_{[N]}$, this becomes \[\abs{\frac{1}{N}\sum_{n=1}^{N}w(n)T^nf}^2 = \frac{1}{N}\left(\frac{1}{N}\sum_{n=1}^{N}\Delta_{0}w(n)T^n\abs{f}^2\right) + \frac{2}{N}\sum_{h=1}^{N-1}\text{Re}\left(\frac{1}{N}\sum_{n=1}^{N}\Delta_{h}w(n)T^n(f\cdot \overline{T^hf})\right),\] and hence by induction for the functions $f\cdot \overline{T^hf} \in L^{2^s}(\mu)$ and weights $\Delta_hw$ we have
    \begin{align}
        \norm{\frac{1}{N}\sum_{n=1}^{N}w(n)T^nf}_{L^{2^{s+1}}(\mu)}^{2} &= \norm{\abs{\frac{1}{N}\sum_{n=1}^{N}w(n)T^nf}^2}_{L^{2^s}(\mu)}\\
        &\lesssim_s \frac{1}{N}\norm{\Delta_{0}w}_{U^s[N]}\norm{\abs{f}^2}_{L^{2^s}(\mu)} + \frac{2}{N} \sum_{h=1}^{N-1}\norm{\Delta_{h}w}_{U^s[N]} \norm{f \cdot\overline{T^hf}}_{L^{2^s}(\mu)}\\
        &\leq \norm{f}_{L^{2^{s+1}}(\mu)}^{2}\left(\frac{1}{N}\norm{\Delta_{0}w}_{U^s[N]} + \frac{2}{N}\sum_{h=1}^{N-1}\norm{\Delta_{h}w}_{U^s[N]}\right),
    \end{align}
    where we use H\"older's Inequality at the end. Recalling that $w=w\one_{[N]}$ we have that \[\norm{\Delta_{h}w}_{U^s[N]} = \norm{\Delta_{-h}w}_{U^s[N]},\] and thus we can reassemble \[\frac{1}{N}\norm{\Delta_{0}w}_{U^s[N]} + \frac{2}{N}\sum_{h=1}^{N-1}\norm{\Delta_{h}w}_{U^s[N]}= \frac{1}{N}\sum_{\abs{h}<N}\norm{\Delta_{h}w}_{U^s[N]}.\] Applying H\"older once again and completing the sum we get \[\frac{1}{N}\sum_{\abs{h}<N}\norm{\Delta_{h}w}_{U^s[N]} \lesssim \left(\frac{1}{N}\sum_{h\in \Z}\norm{\Delta_{h}w}_{U^s[N]}^{2^s}\right)^{1/2^s} = \norm{w}_{U^{s+1}[N]}^{2},\] so we deduce \[\norm{\frac{1}{N}\sum_{n=1}^{N}w(n)T^nf}_{L^{2^{s+1}}(\mu)}^{2}\lesssim_{s+1} \norm{w}_{U^{s+1}[N]}^{2} \norm{f}_{L^{2^{s+1}}(\mu)}^{2},\] which establishes the result.
\end{proof}

\begin{remark}
    Iterating the usual ergodic induction scheme leads to the similar estimate 
    \begin{equation}
    \label{e:Us+1estimate}
        \norm{\frac{1}{N}\sum_{n=1}^{N}w(n)T^nf}_{L^{2^s}(\mu)} \lesssim_{s} \norm{w}_{U^{s+1}[N]}\norm{f}_{L^{2^s}(\mu)},
    \end{equation}
    which does however require a whole additional order of uniformity. The issue lies at the base of the iteration, where van der Corput only gives \[\norm{\frac{1}{N}\sum_{n=1}^{N}w(n)T^nf}_{L^{2}(\mu)} \lesssim \norm{w}_{U^{2}[N]}\norm{f}_{L^{2}(\mu)},\] whereas we start from the sharper \[\norm{\frac{1}{N}\sum_{n=1}^{N}w(n)T^nf}_{L^{4}(\mu)} \lesssim \norm{w}_{U^{2}[N]}\norm{f}_{L^{4}(\mu)}.\]
\end{remark}

By standard methods as in \cite{EKLR,Fan} we get \Cref{p:PositiveResult}. We again include full details.

\begin{proof}[Proof of \Cref{p:PositiveResult}]
    For any fixed $\rho>1$ and lacunary time of the form $\lfloor \rho^m \rfloor$ with $m\geq 0$ \Cref{l:L2SUS} yields \[\norm{\frac{1}{\lfloor \rho^m \rfloor}\sum_{n=1}^{\lfloor \rho^m \rfloor}w(n)T^nf}_{L^{2^{s}}(\mu)}^{2^s} \lesssim_{\rho} \frac{1}{m^{\eta2^s}},\] when $m$ is sufficiently large. Since $\eta>\frac{1}{2^s}$, this is summable, so by Fubini-Tonelli
    \[\int \sum_{m\geq 0}\abs{\frac{1}{\lfloor \rho^m \rfloor}\sum_{n=1}^{\lfloor \rho^m \rfloor}w(n)T^nf}^{2^s} \dd\mu < \infty,\] which implies that \[\frac{1}{\lfloor \rho^{m} \rfloor}\sum_{n=1}^{\lfloor \rho^{m} \rfloor}w(n)f(T^nx) \to 0\] pointwise almost everywhere as $m\to \infty$. For a general time $N\geq 1$ we can find some $m\geq 0$ so that $\lfloor \rho^m \rfloor \leq N < \lfloor\rho^{m+1} \rfloor$ and pointwise bound \[\norm{\frac
    {1}{N}\sum_{n=1}^{N}w(n)f(T^nx) -  \frac{1}{\lfloor \rho^{m} \rfloor}\sum_{n=1}^{\lfloor \rho^{m} \rfloor}w(n)f(T^nx)}_{L^{\infty}(\mu)} \lesssim \norm{w}_{\ell^{\infty}(\Z)}\norm{f}_{L^{\infty}}\cdot \frac{\lfloor\rho^{m+1}\rfloor-\lfloor\rho^m\rfloor}{\lfloor\rho^m\rfloor},\]
    so after sending $N\to \infty$, and hence to $m\to \infty$, \[\limsup_{N\to\infty} \abs{\frac
    {1}{N}\sum_{n=1}^{N}w(n)f(T^nx)} \lesssim \norm{w}_{\ell^{\infty}}\norm{f}_{L^{\infty}(\mu)}(\rho-1).\] Taking $\rho\searrow 1$ in a countable way and consolidating full measure sets finishes the job.
\end{proof}

\begin{remark}
\label{remark:summability}
    More generally, the result follows whenever the summability condition \[\sum_{N=1}^{\infty}\frac{\norm{w}_{U^s[N]}^{2^s}}{N} < \infty,\] holds for some $s\geq 2$, see also \cite[Theorem~2]{Fan}. 
\end{remark}

\subsection{Multiplicative Characters}
\label{subsec:Multiplicative}
Given an odd prime $p\geq 3$, we consider characters $\chi$ of the multiplicative group $\F_{p}^{\times}$, this is, maps from $\F_{p}^{\times}$ to the unit circle such that $\chi(rt)=\chi(r)\chi(t)$ for every $r,t\in \F_{p}^{\times}$. We extend any such map to a map $\F_{p}\to \S \cup \set{0}$ by simply putting $\chi(0)=0$ and call it a \textit{multiplicative character} of $\F_p$. 

The order $q$ of a multiplicative character is the smallest natural number such that $\chi^{q} = 1$ on $\F_{p}^{\times}$. By Lagrange's Theorem, for any multiplicative character $\chi$, we have \[1 \leq q \leq p-1, \qquad q \mid p-1,\] with $q=1$ if and only if $\chi$ is the trivial character $\chi = \one_{\F_{p}^{\times}}$. Conversely, for any $q>1$ as above there exists a character $\chi$ of order $q$ obtained by fixing some $g\in \F_{p}^{\times}$ so that $\F_p^{\times} = \set{g^{j} \colon 1\leq j \leq p-1}$ and putting \[\chi(g^j) = e(j/q),\ 1\leq j \leq p-1.\] For $q=2$, we let $\chi_{p}$ denote the unique multiplicative character of order $q$, which we call the \textit{quadratic character} of $\F_p$ and is given explicitly by \[\chi_p(t) = \begin{cases}
    0,\ &\text{ if }t=0,\\
    1,\ &\text{ if }t\neq 0 \text{ is a square in }\F_p,\\
    -1,\ &\text{ if }t\neq 0 \text{ is not a square in }\F_p.
\end{cases}\]

Some standard facts about non-trivial multiplicative characters, are that they satisfy the following perfect cancellation identity
\begin{equation}
\label{e:cancellationidentity}
    \sum_{t\in \F_{p}}\chi(t) = 0,
\end{equation}
along with the following estimates. The first one is elementary, while the second one is a classical consequence of Weil's bound (see \cite[Corollary~11.24]{IwaniecKowalski}).

\begin{lemma}[Gauss-sum Estimate]
\label{l:GaussSum}
    For every non-trivial multiplicative character $\chi$ and every $\xi\in \F_p$, \[\abs{\sum_{t\in \F_p} \chi(t)e\left(\frac{-\xi t}{p}\right)} \leq p^{1/2}.\]
\end{lemma}

\begin{lemma}[Weil's bound]
\label{l:WeilBound}
    Let $\chi$ be a non-trivial multiplicative character, let $k\geq 1$ and let $s_1,\ldots,s_{2k} \in \F_p$ be pairwise distinct shifts. Then \[\abs{\sum_{t\in \F_p}\chi\left(\prod_{j=1}^{k}(t+s_j)\right)\roverline{\chi\left(\prod_{j=k+1}^{2k}(t+s_j)\right)}} \leq 2kp^{1/2}.\]
\end{lemma}

\section{Character Words}
\label{sec:CharacterWord}

We now introduce a convenient family of objects which allows us to pass the cancellation estimates of the quadratic characters from $\F_p$ to $\Z$. Instead of the usual quadratic character, we  consider the modified quadratic character \[\tau_p(t) \coloneqq \chi_p(t) + \one_{\set{0}}(t),\] and note that \Cref{l:GaussSum} immediately implies \begin{equation}
\label{e:GaussSum}
    \abs{\frac{1}{p}\sum_{t\in \F_p}\tau_{p}(t)e\left(\frac{-\xi t}{p}\right)} \leq \frac{2}{p^{1/2}}
\end{equation}
for every $\xi \in \F_p$.\\

Given $d\geq 1$ a prime $p\geq 3$ and a word $v$ of length $d$ on the symbols $\set{-1,1}$, which we regard as a function $v\colon \Z_{d} \to \set{-1,1}$, we put $\ell=dp$ and define the \textit{character word} $V$ of length $\ell$ induced by $v$, which we regard as a function \[V\colon \Z_{\ell}\to \set{-1,1},\] via \[V(dt+r) = \tau_{p}(t)v(r),\] for $0\leq t <p$ and $0\leq r < d$. Informally, we are decomposing our ambient space into a tensor of arithmetic progression/cyclic subgroup; this is precisely the block product $v\times \tau_{p}$ used in \cite{Keane} with $\tau_p$ prescribing the signs of the concatenated copies of $v$. 

\begin{ex}
    If $d=4$, $p=3$ and \[v = (1,-1,-1,1),\] then since $\tau_3(0)=\tau_{3}(1)=1$ and $\tau_3(2)=-1$, we obtain the character word \[V = \bigl( \underbrace{1,-1,-1,1}_{t=0}, \underbrace{1,-1,-1,1}_{t=1}, \underbrace{-1,1,1,-1}_{t=2}\bigr).\]
\end{ex}

\subsection{Fourier Cancellation}

Character words retain the Fourier cancellation of quadratic characters; here and throughout we let \[\widehat{V}\colon \Z_{\ell}\to \C, \qquad \widehat{V}(\xi ) \coloneqq \frac{1}{\ell}\sum_{n=0}^{\ell-1}V(n)e\left(\frac{-\xi n}{\ell}\right)\] denote the Fourier Transform of $V\colon \Z_{\ell} \to \C$.

\begin{lemma}
\label{l:FourierCharacterWord}
    For any $d\geq 1$, any prime $p\geq 3$ and any word $v\colon \Z_{d}\to \set{-1,1}$ the character word $V$ of length $\ell=dp$ induced by $v$ satisfies \[\max_{\xi \in \Z_{\ell}} \abs{\widehat{V}(\xi )} \leq \frac{2}{p^{1/2}}.\]
\end{lemma}

\begin{proof}
    Decomposing each $n$ as $n=dt+r$, with $0\leq t <p$ and $0\leq r < d$ and recalling that $\ell=dp$, we may split the sum as follows
    \begin{align}
        \widehat{V}(\xi ) = \frac{1}{\ell}\sum_{n=0}^{\ell-1}V(n)e\left(\frac{-\xi n}{\ell}\right) &= \frac{1}{dp}\sum_{t=0}^{p-1}\sum_{r=0}^{d-1}\tau_{p}(t)v(r)e\left(\frac{-\xi (dt+r)}{dp}\right)\\
        &= \left(\frac{1}{p}\sum_{t=0}^{p-1}\tau_p(t)e\left(\frac{-\xi t}{p}\right)\right)\left(\frac{1}{d}\sum_{r=0}^{d-1}v(r)e\left(\frac{-\xi r}{dp}\right)\right).
    \end{align}
    After taking absolute values, the first term is bounded by the modified Gauss-sum estimate \eqref{e:GaussSum}, and the second one is trivially bounded by $1$, giving the result.
\end{proof}

In particular, by summation by parts we can pass the estimate to arbitrary intervals, which is what we will crucially use later on.

\begin{lemma}
\label{l:FourierDecay}
    Let $\ell \geq 1$ and let $V\colon \Z\to \C$ be an $\ell$-periodic function such that \[\max_{\xi \in \Z_{\ell}}\abs{\widehat{V}(\xi )} \leq \delta.\] Then, for every interval $I\subseteq \Z$, \[\sup_{\beta\in \T}\abs{\sum_{n\in I}V(n)e(n\beta)} \lesssim \delta(\abs{I}+\ell\log \ell).\] 
\end{lemma}

\begin{proof}
    By Fourier Inversion, for $0\leq n < \ell$ we may write \[V(n)=\sum_{\xi =0}^{\ell-1}\widehat{V}(\xi )e\left(\frac{\xi n}{\ell}\right)\] so the trivial bound for exponential sum implies
    \begin{align}
        \abs{\sum_{n\in I}V(n)e(n\beta)} \leq \delta \sum_{\xi =0}^{\ell-1}\min\set{\abs{I}, \frac{1}{2\norm{\beta+\xi /\ell}_{\T}}}.
    \end{align}
    Since the collection of frequencies \[\set{\xi /\ell \colon 0 \leq \xi  < \ell}\] is $1/\ell$-separated on $\T$, for every $\beta \in \T$ and every $0\leq j \leq \ell-1$ there are at most two values of $\xi $ so that \[\frac{j}{2\ell}\leq\norm{\beta + \xi /\ell}_{\T} < \frac{j+1}{2\ell}.\] Splitting the sum accordingly and using the first bound for the $j=0$ term and the second bound for the remaining ones, we get
    \begin{align}
        \abs{\sum_{n\in I}V(n)e(n\beta)} \lesssim \delta\abs{I} + \delta\sum_{j=1}^{\ell-1}\frac{2\ell}{2j} \lesssim \delta(\abs{I} + \ell \log \ell).
    \end{align}
\end{proof}

\subsection{Gowers Norm Estimates} In the same line as above, uniformity properties of the quadratic characters can be passed on to $\Z_{\ell}$.

\begin{proposition}
\label{p:GowersNormCharacterWord}
    For any $d\geq 1$, any prime $p\geq 3$ and any word $v\colon \Z_{d}\to \set{-1,1}$ the character word $V$ of length $\ell=dp$ induced by $v$ satisfies \[\norm{V}_{U^s(\Z_{\ell})} \lesssim_s \frac{1}{p^{1/2^s}},\] for every $s\geq 2$.
\end{proposition}

To establish this proposition we introduce a further character word $V_{\chi}\colon \Z_{\ell}\to \set{-1,0,1}$ defined using the actual quadratic character as \[V_{\chi}(dt+r)=\chi_p(t)v(r),\] for $0\leq t<p$, $0\leq r <d$. By the triangle inequality for Gowers norms we then have \[\norm{V}_{U^s(\Z_{\ell})} \leq \norm{V-V_{\chi}}_{U^s(\Z_{\ell})} +\norm{V_{\chi}}_{U^s(\Z_{\ell})} ,\] so \Cref{p:GowersNormCharacterWord} is reduced to the following two lemmas.

\begin{lemma}\label{l:CharacterWordGowers0} 
    For any $d\geq 1$, any prime $p\geq 3$ and any word $v\colon \Z_{d}\to \set{-1,1}$, the character words $V$ and $V_{\chi}$ of length $\ell=dp$ induced by $v$ satisfy \[\norm{V-V_{\chi}}_{U^s(\Z_{\ell})} \leq \frac{1}{p^{1/2^s}},\] for every $s\geq 2$.
\end{lemma}

\begin{lemma}
\label{l:CharacterWordGowers}
    For any $d\geq 1$, any prime $p\geq 3$ and any word $v\colon \Z_{d}\to \set{-1,1}$, the character word $V_{\chi}$ of length $\ell=dp$ induced by $v$ satisfies \[\norm{V_{\chi}}_{U^s(\Z_{\ell})} \lesssim_s \frac{1}{p^{1/2^s}},\] for every $s\geq 2$.
\end{lemma}

We address them in order.

\begin{proof}[Proof of \Cref{l:CharacterWordGowers0}]
    We simply exploit the identity \[\abs{(V-V_{\chi})(n)} = \one_{\set{0 \leq n \ (\text{mod } \ell) \leq d-1}}\] which after expanding the Gowers norm reduces matters to a counting problem
    \begin{align}
        &\norm{V-V_{\chi}}_{U^s(\Z_{\ell})}^{2^s}\\
        &\qquad \leq \frac{1}{\ell^{s+1}} \abs{\set{(n,h_1,\ldots,h_s) \colon  0\leq n+\omega \cdot h \ (\text{mod } \ell) \leq d-1,\ \forall \omega\in\set{0,1}^{s}}} \\
        &\qquad \qquad \leq \frac{1}{\ell^{s+1}} \abs{\set{(n,h_1,\ldots,h_s) \colon  0\leq n \ (\text{mod } \ell) \leq d-1}},
    \end{align}
    which yields the bound
    \[\norm{V-V_{\chi}}_{U^s(\Z_{\ell})}^{2^s} \leq \frac{1}{\ell^{s+1}} \cdot d \cdot \ell^{s} = \frac{1}{p},\] since $\ell = dp$.
\end{proof}

\begin{proof}[Proof of \Cref{l:CharacterWordGowers}]
    Begin by expanding \[\norm{V_{\chi}}_{U^s(\Z_{\ell})}^{2^s} = \frac{1}{\ell^{s+1}}\sum_{n,h_1,\ldots,h_s \in \Z_{\ell}} \prod_{\omega \in \set{0,1}^{s}}\CC^{\abs{\omega}}V_{\chi}(n+\omega \cdot h).\] Since $\ell=dp$, after naturally injecting $\F_p$ into $\Z_{\ell}$, we may uniquely write  $n=dt+r$ and $h_i = dt_i+r_i$, with $t,t_1,\ldots,t_s\in \F_p$ and $0\leq r,r_1,\ldots,r_s < d$ and split \[\norm{V_{\chi}}_{U^s(\Z_{\ell})}^{2^s} = \frac{1}{p^{s+1}}\frac{1}{d^{s+1}}\sum_{t,t_1,\ldots,t_s \in \F_p}\sum_{r,r_1,\ldots,r_s = 0}^{d-1}\prod_{\omega\in \set{0,1}^{s}}V_{\chi}\left(dt+r + \sum_{i=1}^{s}\omega_i(dt_i+r_i)\right).\] To apply the definition of $V_{\chi}$ we now need to fully decompose \[n + \omega \cdot h = dt + r + \sum_{i=1}^{s}\omega_i(dt_i+r_i) = d(t+\omega \cdot \Vec{t}) + (r+\omega \cdot \Vec{r}),\] for which it remains to reduce $r+\omega\cdot \Vec{r}$ modulo $d$. Letting \[t_{\omega} \coloneqq \left\lfloor \frac{r+\omega\cdot \Vec{r}}{d} \right\rfloor, \quad r_{\omega} \coloneqq r+\omega \cdot \Vec{r} - dt_{\omega},\] we have \[n+\omega \cdot h = d(t+\omega \cdot \Vec{t} + t_{\omega}) + r_{\omega} \mod \ell\] with $0\leq r_{\omega} < d$. Regarding the quantity  $t+\omega \cdot \Vec{t} + t_{\omega}$ modulo $p$, we may now unambiguously write 
    \begin{align}
        \norm{V_{\chi}}_{U^s(\Z_{\ell})}^{2^s} &= \frac{1}{p^{s+1}}\frac{1}{d^{s+1}}\sum_{t,t_1,\ldots,t_s\in \F_p}\sum_{r,r_1,\ldots,r_s = 0}^{d-1}\prod_{\omega\in \set{0,1}^{s}}\chi_p(t+\omega\cdot \Vec{t}+t_{\omega})v(r_{\omega})\\
        &=\frac{1}{d^{s+1}}\sum_{r,r_1,\ldots,r_s=0}^{d-1}\left(\prod_{\omega\in \set{0,1}^{s}}v(r_{\omega})\right) A(r,r_1,\ldots,r_s),\label{e:V0norm}
    \end{align}
    where we note that $r_{\omega}$ only depends on $r,r_1,\ldots,r_s$ and put \[A = A(r,r_1,\ldots,r_s) \coloneqq \frac{1}{p^{s+1}}\sum_{t,t_1,\ldots,t_s \in \F_p}\prod_{\omega\in \set{0,1}^{s}}\chi_p(t+\omega\cdot \Vec{t}+t_{\omega}).\] 
    Since $|v| \leq 1$, matters reduce to bounding
    \begin{align}
        |A(r,r_1,\dots,r_s)| \lesssim_s \frac{1}{p}
    \end{align}
    uniformly for $0 \leq r,r_1,\dots,r_s < d$. To do so, we write $\omega=(\eta,i)$ for $\eta\in \set{0,1}^{s-1}$ and $i\in \set{0,1}$, so that we may express \[\prod_{\omega \in \set{0,1}^{s}}\chi_p(t+\omega \cdot \Vec{t} + t_{\omega}) = \prod_{\eta\in \set{0,1}^{s-1}}\chi_p(t+\eta \cdot \Vec{t'} + t_{(\eta,0)})\chi_p(t+\eta \cdot \Vec{t'} + t_{s} + t_{(\eta,1)}),\] where $\Vec{t'}=(t_1,\ldots,t_{s-1})$. After averaging over the last cube direction $t_s$ with the change of variables $(t,t_s)\mapsto (t,t+t_s)$, this allows us to write \[A = \frac{1}{p^{s-1}}\sum_{t_1,\ldots,t_{s-1}\in \F_p}A_0(t_1,\ldots,t_{s-1}) A_1(t_1,\ldots,t_{s-1}),\] where
    \begin{align}
        A_0(t_1,\ldots,t_{s-1}) &\coloneqq \frac{1}{p}\sum_{t\in \F_p}\prod_{\eta\in \set{0,1}^{s-1}}\chi_p(t+\eta \cdot \Vec{t'} + t_{(\eta,0)}),\\
        A_1(t_1,\ldots,t_{s-1}) &\coloneqq \frac{1}{p}\sum_{t_s\in \F_p}\prod_{\eta\in \set{0,1}^{s-1}}\chi_p(t_s+\eta \cdot \Vec{t'} + t_{(\eta,1)}),
    \end{align}
    We now estimate these quantities using \Cref{l:WeilBound}. First, if $t_1,\ldots,t_{s-1}$ are so that all shifts $\eta \cdot \Vec{t'} + t_{(\eta,0)}$ are distinct modulo $p$, then \[\abs{A_0(t_1,\ldots,t_{s-1})} \lesssim_s \frac{1}{p^{1/2}},\] and if additionally $\eta \cdot \Vec{t'} + t_{(\eta,1)}$ are also all distinct modulo $p$, then similarly \[\abs{A_1(t_1,\ldots,t_{s-1})} \lesssim_s \frac{1}{p^{1/2}}.\] It follows that the contribution of these $t_1,\ldots,t_{s-1}$ to the average $A$ is $\lesssim_{s} 1/p$. For the remaining $t_1,\ldots,t_{s-1}$, bound $A_0$ and $A_{1}$ by $1$ and observe that for both $i\in \set{0,1}$ there are at most $\lesssim_{s} p^{s-2}$ possible $t_1,\ldots,t_{s-1}$ so that one of the shifts $\eta \cdot \Vec{t'} + t_{(\eta,i)}$ repeats. This means that we can bound their contribution to the average $A$ by $1/p$ and get \[\abs{A(r,r_1,\ldots,r_s)} \lesssim_{s} \frac{1}{p},\]
    as desired.
\end{proof}

\Cref{p:GowersNormCharacterWord} can be transferred to arbitrary intervals in $\Z$ thanks to the following result.

\begin{lemma}
\label{l:GowersNormPeriodic}
    For any $s\geq 2$, $\ell\geq 2$, any interval $I\subseteq \Z$ and any $1$-bounded, $\ell$-periodic $V\colon \Z\to \C$ one has \[\norm{V}_{U^s[I]}^{2^s} \lesssim_s \norm{V}_{U^s(\Z_{\ell})}^{2^s} + \frac{\ell}{\abs{I}}.\] 
\end{lemma}

\begin{proof}
    Since the cyclic Gowers norms are translation invariant, we may translate $I$ and assume without loss of generality that $I = [M]$. Then, by definition \[\norm{V}_{U^{s}(\Z_{\ell})}^{2^s} = \frac{1}{\ell^{s+1}}\sum_{n,h_1,\ldots,h_s\in \Z_{\ell}}\prod_{\omega\in \set{0,1}^{s}}\CC^{\abs{\omega}}V(n+\omega\cdot h),\] meanwhile \[\norm{V}_{U^s[M]}^{2^s} = \frac{1}{M^{s+1}}\sum_{\substack{n,h_1,\ldots,h_s \in \Z \\ n+\omega\cdot h \in [M],\ \forall \omega \in \set{0,1}^{s}}}\prod_{\omega\in \set{0,1}^{s}}\CC^{\abs{\omega}}V(n+\omega\cdot h).\] For simplicity we let \[\II \coloneqq \set{(n,h_1,\ldots,h_s)\in \Z^{s+1} \colon n+\omega\cdot h \in [M],\ \forall \omega\in \set{0,1}^{s}}.\] The key observation is that since $V$ is $\ell$-periodic, the map \[(n,h_1,\ldots,h_s) \mapsto V(n+\omega\cdot h)\] is $\ell$-periodic in every entry, so for any block of the form 
    \begin{equation}
    \label{e:blocks}
        B = B(z) = \set{0,\ldots,\ell-1}^{s+1} + \ell z,\ z\in \Z^{s+1}
    \end{equation}
    one has 
    \begin{equation}
    \label{e:blocknorm}
        \norm{V}_{U^{s}(\Z_{\ell})}^{2^s} = \frac{1}{\ell^{s+1}}\sum_{(n,h_1,\ldots,h_s)\in B} \prod_{\omega\in \set{0,1}^{s}}\CC^{\abs{\omega}}V(n+\omega\cdot h).
    \end{equation}
    Covering $\II$ by blocks of the form \eqref{e:blocks}, we identify two types of blocks, those that are entirely contained in $\II$ and those that protrude outside of $\II$; we denote the set of tuples $(n,h_1,\ldots,h_s) \in \II$ that belong to a block $B$ with $B\not\subseteq \II$ by $\RR$, and note that $\RR$ is contained inside \[\set{(n,h_1,\ldots,h_s)\in \mathcal{I} \colon \min_{\omega\in \set{0,1}^{s}} \dist(n+\omega \cdot h, \set{0,M}) \leq (s+1)\ell}\]
    which has size $\lesssim_s \ell M^{s}$. If we now split 
    \begin{multline}
        \norm{V}_{U^s[M]}^{2^s} = \frac{1}{M^{s+1}}\sum_{(n,h_1,\ldots,h_s)\in \II\setminus \RR}\prod_{\omega\in \set{0,1}^{s}}\CC^{\abs{\omega}}V(n+\omega\cdot h)\\
        + \frac{1}{M^{s+1}}\sum_{(n,h_1,\ldots,h_s)\in \RR}\prod_{\omega\in \set{0,1}^{s}}\CC^{\abs{\omega}}V(n+\omega\cdot h),
    \end{multline}
we may bound the second term by $\lesssim_s \frac{\ell}{M}$, while the first term can be bounded as \[\frac{1}{M^{s+1}}\sum_{(n,h_1,\ldots,h_s)\in \II\setminus \RR}\prod_{\omega\in \set{0,1}^{s}}\CC^{\abs{\omega}}V(n+\omega\cdot h) \lesssim \norm{V}_{U^{s}(\Z_{\ell})}^{2^s}.\]
\end{proof}

Combining \Cref{p:GowersNormCharacterWord} and \Cref{l:GowersNormPeriodic} immediately produces the following corollary.

\begin{cor}
\label{c:GowersNormsCharacterWord}
    Let $s\geq 2$, $d\geq 1$, $p\geq 3$ be a prime, $v\colon \Z_d\to \set{-1,1}$ and $V$ be the character word of length $\ell=dp$ induced by $v$, extended $\ell$-periodically to all of $\Z$. Then for any interval $I\subseteq \Z$ with $\abs{I} \geq \ell^2$, \[\norm{V}_{U^s[I]}^{2^s} \lesssim_s \frac{1}{p}.\] Moreover, if $I\subseteq [N]$, then \[\norm{V\cdot \one_{I}}_{U^s[N]} \lesssim_{s} \frac{1}{p^{1/2^s}}\left(\frac{\abs{I}}{N}\right)^{\frac{s+1}{2^s}}.\]
\end{cor}

\begin{remark}\label{rem:symmetry}
    In the next section we will not only work with character words but with their signed translates $V_{\eps,a}$ defined for each sign $\eps\in\set{-1,1}$ and each shift $a\in \Z$ via \[V_{\eps,a}(n) = \eps V(n+a).\] Since the estimates above rely at their core on cyclic Gowers norms, which are invariant under translations and sign changes, both hold uniformly in signs and shifts.
\end{remark}

\section{Construction of the Weight}
\label{sec:Construction}

\subsection{Construction}

Throughout, we fix a non-decreasing function $\psi\colon \R_{+}\to \R_{+}$ with $\psi(x) \to \infty$ as $x\to \infty$. We construct a sequence of (rapidly) increasing \say{master intervals}, 
\[ I^{(k)} := (A_k,B_k] \cap \Z,\  k \geq 1  \]
inside of which we embed many copies of nested character words. The construction is inductive and involved, so we begin with an informal description:

Suppose we have constructed a character word up to time $k$, call it $V_k$, whose length is $\ell_k$. To extend the construction to time $k+1$, we start by repeating the word $V_k$ $q_{k+1}$ many times, where we think of $q_{k+1}$ as being enormous relative to all earlier parameters, to arrive at an intermediate word
\[ v_{k+1}(n) := \sum_{0 \leq j < q_{k+1}} V_k(n - j \ell_k) \cdot \one_{[j\ell_k,(j+1)\ell_k)}(n); \]
we set $d_{k+1} := \ell_k q_{k+1}$. We now apply our tensor perspective to lift our intermediate word $v_{k+1}$ to the $k+1$-th level character word $V_{k+1}$: with $p_{k+1}$ a sufficiently large prime (certainly larger than $ k^{100}$), we define
\[ V_{k+1}(d_{k+1} t + r) := \tau_{p_{k+1}}(t) v_{k+1}(r), \qquad 0 \leq t < p_{k+1}, \ 0 \leq r < d_{k+1}. \]
The following schematic illustrates this process.

\begin{center}
\normalsize
\setlength{\unitlength}{1pt}
\makebox[\linewidth][c]{%
\begin{picture}(390,168)
  % Repeated copies of V_k.
  \put(216,161){\makebox(0,0){$q_{k+1}$ copies,\qquad
      $d_{k+1}=q_{k+1}\ell_k$}}
  \put(101,151){\line(1,0){230}}
  \put(101,147){\line(0,1){4}}
  \put(331,147){\line(0,1){4}}
  \put(91,130.5){\makebox(0,0)[r]{$v_{k+1}=$}}
  {\linethickness{0.55pt}%
   \put(101,118){\framebox(230,25){}}
   \multiput(147,118)(46,0){4}{\line(0,1){25}}}%
  \multiput(124,130.5)(46,0){3}{\makebox(0,0){$V_k$}}
  \put(262,130.5){\makebox(0,0){$\cdots$}}
  \put(308,130.5){\makebox(0,0){$V_k$}}

  % Entrywise multiplication by the sign of the t-th copy.
  {\color{teal!60}%
   \put(101,118){\line(1,-1){65}}
   \put(331,118){\line(-1,-1){65}}}%
  \put(216,101){\makebox(0,0){multiply every entry by
      $\tau_{p_{k+1}}(t)$.}}
  {\color{teal}%
   \put(216,85){\vector(0,-1){24}}}%

  % Signed copies of v_{k+1}. Each labelled cell is 100pt wide.
  \put(37,40.5){\makebox(0,0)[r]{$V_{k+1}=$}}
  \put(166,28){\color{teal!10}\rule{100\unitlength}{25\unitlength}}
  {\linethickness{0.55pt}%
   \put(43,28){\line(1,0){346}}
   \put(43,53){\line(1,0){346}}
   \put(43,28){\line(0,1){25}}
   \put(143,28){\line(0,1){25}}
   \put(289,28){\line(0,1){25}}
   \put(389,28){\line(0,1){25}}
   {\color{teal}%
    \put(166,28){\line(1,0){100}}
    \put(166,53){\line(1,0){100}}
    \put(166,28){\line(0,1){25}}
    \put(266,28){\line(0,1){25}}}}%
  \put(93,40.5){\makebox(0,0){$\tau_{p_{k+1}}(0)\,v_{k+1}$}}
  \put(154.5,40.5){\makebox(0,0){$\cdots$}}
  \put(216,40.5){\makebox(0,0){$\tau_{p_{k+1}}(t)\,v_{k+1}$}}
  \put(277.5,40.5){\makebox(0,0){$\cdots$}}
  \put(339,40.5){\makebox(0,0){$\tau_{p_{k+1}}(p_{k+1}{-}1)v_{k+1}$}}
  \put(43,20){\line(1,0){346}}
  \put(43,20){\line(0,1){4}}
  \put(389,20){\line(0,1){4}}
  \put(216,9){\makebox(0,0){$p_{k+1}$ signed copies,\qquad
      $\ell_{k+1}=p_{k+1}d_{k+1}$}}
\end{picture}%
}
\end{center}

Further ahead, we will concatenate character words throughout the master intervals $I^{(k)}$ to eventually force non-convergence phenomena to repeat sufficiently often.
\medskip

We now describe the construction rigorously:\\

Initialize $d_0=1$, $\ell_0=1$, $B_0=1$ and $V_0=1$ regarded as a word of length $\ell_0$ on the symbols $\set{-1,1}$. At stage $k\geq 1$, given $\ell_{k-1}$, $B_{k-1}$ and $V_{k-1}$, put
\begin{equation}
\label{e:defi-qkdk}
    q_k \coloneqq \left \lceil\frac{2^kB_{k-1}}{\ell_{k-1}} \right \rceil, \qquad d_k \coloneqq q_k\ell_{k-1}, 
\end{equation}
and define $v_k$ to be the word of length $d_k$ consisting of $q_k$ copies of $V_{k-1}$, which itself is a word of length $\ell_{k-1}$. Choose an odd prime $p_k$ sufficiently large so that 
\begin{equation}
\label{e:defi-pk}
    d_k(1+\log(k+1)) \leq \psi(\log(d_k p_k))
\end{equation}
and set
\begin{equation}
\label{e:defi-lkAkBk}
    \ell_k \coloneqq d_kp_k, \qquad A_{k}\coloneqq \ell_k^{4}, \qquad B_k = A_k(k+1)^{2\ell_k}.
\end{equation}
Finally, let $V_k$ be the character word of length $\ell_k$ induced by $v_k$, this is,
\begin{equation}
\label{e:defi-Vk}
    V_{k}(d_kt + r) = \tau_{p_k}(t)v_{k}(r),
\end{equation}
for $0\leq t<p_k$ and $0\leq r<d_k$. 
\\

Carrying out this procedure we obtain sequences $(\ell_k)_{k\geq 0}$, $(d_k)_{k\geq 1}$, $(A_k)_{k\geq 1}$ and $(B_{k})_{k\geq 0}$, such that for every $k\geq 1$
\begin{equation}
    d_k \geq 2^kB_{k-1}, \qquad \ell_k \geq 3 \cdot 2^{k}B_{k-1},
\end{equation}
and also
\begin{equation}
    B_{k-1} \leq d_k  < \ell_k < A_k < B_k, \qquad d_{k-1} \mid \ell_{k-1} \mid d_k.
\end{equation}

Regarding each word $V_k$ as a function on $\Z$ by considering its $\ell_k$-periodic extension  we crucially have the following nesting property, established by induction (see also the schematic).

\begin{lemma}
\label{l:Restriction}
    Let $k\geq 0$ and $j\geq k+1$. Then, for every $m\in \Z$ there exists a sign $\eps_{j,m}\in \set{-1,1}$, such that \[V_j(n)=\eps_{j,m}V_k(n),\] for $md_{k+1} \leq n < (m+1)d_{k+1}$.
\end{lemma}

\begin{proof}
    For $j=k+1$, let $m\in \Z$ and write each $md_{k+1} \leq n < (m+1)d_{k+1}$ as \[n=md_{k+1}+r\] with $0\leq r < d_{k+1}$. Decomposing $m=ap_{k+1}+t$ with $a\in \Z$ and $0\leq t<p_{k+1}$, we have that \[n = ap_{k+1}d_{k+1} + td_{k+1} + r = a\ell_{k+1} + td_{k+1} + r,\] and since $V_{k+1}$ is $\ell_{k+1}$-periodic \[V_{k+1}(n) = V_{k+1}(td_{k+1} + r).\] Since $0\leq t<p_{k+1}$ and $0\leq r < d_{k+1}$, by definition we now have \[V_{k+1}(n) = \tau_{p_{k+1}}(t)v_{k+1}(r).\] Letting $\eps_{k+1,m} = \tau_{p_{k+1}}(t)$ and recalling that $v_{k+1}$ consists of $q_{k+1}$ copies of $V_k$, we get \[V_{k+1}(n) = \eps_{k+1,m}V_{k}(r).\] Finally, because $n=md_{k+1}+r$ and $d_{k+1}=q_{k+1}\ell_k$, we have that $r \equiv n \mod \ell_k$, so, since $V_k$ is $\ell_k$-periodic \[V_{k+1}(n) = \eps_{k+1,m}V_{k}(n)\] as claimed. For the inductive step, assume we have \[V_{j}(n) = \eps_{j,m}V_{k}(n)\] for $md_{k+1} \leq n < (m+1)d_{k+1}$. Since $d_{k+1} \mid d_{j+1}$, there exists some $a \in \Z$ such that the interval \[[md_{k+1},(m+1)d_{k+1}) \cap \Z\]is entirely contained in the interval \[[ad_{j+1},(a+1)d_{j+1}) \cap \Z.\] Proceeding exactly as above, we can obtain a sign $\eta$, such that \[V_{j+1}(n)=\eta V_j(n),\] for $ad_{j+1}\leq n < (a+1)d_{j+1}$, and hence, putting $\eps_{j+1,m} = \eta \eps_{j,m}$, \[V_{j+1}(n) = \eps_{j+1,m}V_{k}(n),\] for $md_{k+1} \leq n < (m+1)d_{k+1}$.
\end{proof}

We are now ready to define our weight $w$. For each $k\geq 1$, fix an enumeration of the collection of pairs of signs and shifts \begin{align}\label{e:signtrans}
\EE_k \coloneqq \set{-1,1}\times \set{0,\ldots,\ell_k-1} = \set{(\eps^{(k)}_j,a^{(k)}_j) \colon 1\leq j \leq 2\ell_k}
\end{align}
and consider the master interval \[I^{(k)} \coloneqq (A_k,B_k] \cap \Z,\] subdivided into the disjoint intervals \[I^{(k)}_{j} \coloneqq (A_k(k+1)^{j-1}, A_k(k+1)^{j}] \cap \Z, \qquad \big \lvert I^{(k)}_{j} \big \rvert = A_k(k+1)^{j-1}k\] for $1\leq j\leq 2\ell_k$. Finally define $w\colon \Z\to \set{-1,0,1}$ via 
\begin{equation}
\label{e:defiweight}
    w(n) \coloneqq \begin{cases}
    \eps^{(k)}_{j}V_k(n+a^{(k)}_j),\ &\text{ if } n\in I^{(k)}_{j},\ \text{ for some } k\geq 1 \text{ and } 1\leq j\leq 2\ell_{k},\\
    0,\ &\text{ otherwise;}
\end{cases}
\end{equation}
the below schematic describes the behavior of $w$ on its support, $\bigcup_{k \geq 1} I^{(k)}$.

\begin{center}
\small
\setlength{\unitlength}{1pt}
\begin{picture}(350,156)
  \put(175,147){\makebox(0,0){$I^{(k)}$}}
  \put(38,133){\line(1,0){274}}
  \put(38,129){\line(0,1){4}}
  \put(312,129){\line(0,1){4}}

  \put(38,105){\framebox(40,22){$I^{(k)}_1$}}
  \put(78,105){\framebox(50,22){$\cdots$}}
  \put(128,105){\framebox(68,22){$I^{(k)}_j$}}
  \put(196,105){\framebox(32,22){$\cdots$}}
  \put(228,105){\framebox(84,22){$I^{(k)}_{2\ell_k}$}}
  \put(38,94){\makebox(0,0){$A_k$}}
  \put(312,94){\makebox(0,0){$B_k$}}

  % One period inside the j-th interval.
  \put(128,105){\line(-3,-1){78}}
  \put(196,105){\line(3,-1){78}}
  \put(50,55){\framebox(44,24){$\cdots$}}
  \put(230,55){\framebox(44,24){$\cdots$}}
  {\linethickness{0.8pt}
   \put(94,55){\framebox(136,24){$w(n)=\eps^{(k)}_j V_k(n+a^{(k)}_j)$}}}

  \put(94,47){\line(1,0){136}}
  \put(94,47){\line(0,1){4}}
  \put(230,47){\line(0,1){4}}
  \put(162,36){\makebox(0,0){$\ell_k$}}
\end{picture}
\end{center}

The key observation which we will use to demonstrate our divergence phenomenon is as follows:

Suppose we are given a sequence $u = (u_n)_{n\in \Z}\in \set{-1,1}^{\Z}$ so that for infinitely many $k\geq 1$ one has \[u_n=\eps V_k(n+a)\] throughout $[B_k]$ for some $(\eps,a)\in\EE_k$. In this case $w(n)u_n$ is identically $1$ on the block indexed by $(\eps,a)$ and identically $-1$ on the block indexed by $(-\eps,a)$. At the right endpoint $N=A_k(k+1)^j$ of either block $I_{j}^{(k)}$, said block occupies most of the full interval $[N]$ ---indeed, we have $|I_j^{(k)}|/N=k/(k+1)$--- and hence
\begin{align}
    \frac{1}{N} \sum_{n=1}^{N} w(n) u_n = \pm 1 + \epsilon_{u,k}, \qquad |\epsilon_{u,k}| \leq \frac{2}{k+1},
\end{align}
for this $N$; in particular we have a mechanism to arrange for our averages with respect to $u$ to be arbitrarily close to $\pm 1$ infinitely often.

\subsection{Cancellation}

We begin by looking at Fourier cancellation.

\begin{lemma}
\label{l:FourierCancellation}
    For every $k\geq 1$ and $N\in I^{(k)}$ the weight defined by \eqref{e:defiweight} satisfies \[\sup_{\beta\in \T}\abs{\frac{1}{N}\sum_{n=1}^{N}w(n)e(n\beta)} \lesssim \frac{1}{p_k^{1/2}}.\]
\end{lemma}

\begin{proof}
    Let $1\leq J \leq 2\ell_k$ so that $N\in I^{(k)}_{J}$ and split 
    \begin{multline}
        \abs{\frac{1}{N}\sum_{n=1}^{N}w(n)e(n\beta)} \leq \frac{1}{N}\abs{\sum_{n=1}^{B_{k-1}}w(n)e(n\beta)}\\
        + \frac{1}{N}\sum_{j=1}^{J-1}\abs{\sum_{n\in I^{(k)}_{j}}w(n)e(n\beta)} + \frac{1}{N}\abs{\sum_{n\in I^{(k)}_{J} \cap [N]}w(n)e(n\beta)},
    \end{multline}
    using that $w\equiv 0$ between $B_{k-1}$ and $A_k$. For the first contribution, we simply use that $w$ takes values in $\set{-1,0,1}$ and bound \[\frac{1}{N}\abs{\sum_{n=1}^{B_{k-1}}w(n)e(n\beta)} \leq \frac{B_{k-1}}{N}.\] For each full block $I_{j}^{(k)}$, $1\leq j < J$, by definition, \[w(n) = \eps^{(k)}_{j}V_k(n+a^{(k)}_{j}).\] Now, since $V_k$ is $\ell_k$-periodic, so is \[n\mapsto \eps^{(k)}_{j}V_k(n+a^{(k)}_{j}),\] and moreover, \[\widehat{\eps^{(k)}_{j}V_k(\cdot + a^{(k)}_{j})}(b) = \eps^{(k)}_{j}e\left(\frac{a^{(k)}_{j}b}{\ell_k}\right)\widehat{V_k}(b),\] for every $b\in \Z$. By \Cref{l:FourierCharacterWord} we bound \[\max_{b\in \Z} \abs{\widehat{\eps^{(k)}_{j}V_k(\cdot + a^{(k)}_{j})}(b)} = \max_{b\in \Z}\abs{\widehat{V_k}(b)} \leq \frac{2}{p_k^{1/2}},\] and so by \Cref{l:FourierDecay}, \[\sup_{\beta\in \T}\abs{\sum_{n\in I^{(k)}_{j}}w(n)e(n\beta)} \lesssim \frac{1}{p_k^{1/2}}\left(\abs{I_{j}^{(k)}}+\ell_k\log(\ell_k)\right).\] The same argument applies to the partial interval, 
    \[\sup_{\beta\in \T}\abs{\sum_{n\in I^{(k)}_{J} \cap [N]}w(n)e(n\beta)} \lesssim \frac{1}{p_k^{1/2}}\left(\abs{I^{(k)}_{J} \cap [N]} + \ell_k\log(\ell_k)\right),\]
    so putting everything together we obtain the bound \[\abs{\frac{1}{N}\sum_{n=1}^{N}w(n)e(n\beta)} \lesssim \frac{B_{k-1}}{N} + \frac{1}{Np^{1/2}}\left(\sum_{j=1}^{J-1}\left(\abs{I_{j}^{(k)}} + \ell_k\log(\ell_k)\right) + \abs{I^{(k)}_{J} \cap [N]} + \ell_k\log(\ell_k)\right).\] Since the sub-intervals are disjoint we clearly have \[\sum_{j=1}^{J-1}\abs{I_{j}^{(k)}} + \abs{I^{(k)}_{J} \cap [N]} \leq N,\] and since $J \leq 2\ell_k$, \[(J-1)\ell_k\log(\ell_k) + \ell_k\log(\ell_k) \leq \ell_k^{3},\] so after reorganizing we obtain the bound \[\abs{\frac{1}{N}\sum_{n=1}^{N}w(n)e(n\beta)} \lesssim \frac{B_{k-1}}{N} + \frac{1}{Np_k^{1/2}}\left(N + \ell_k^{3}\right).\] We address the first term: since $N\in I^{(k)}$, we have $N\geq A_k = \ell_{k}^{4}$, so, coupled with the fact that $B_{k-1} \leq \ell_k$, we get \[\frac{B_{k-1}}{N} \leq \frac{\ell_{k}}{\ell_{k}^{4}} \leq \frac{1}{\ell_{k}^{3}} \leq \frac{1}{p_{k}^{3}} \leq \frac{1}{p_{k}^{1/2}},\] where we recall that $\ell_k = d_kp_k$. For the second term, \[\frac{1}{N}(N + \ell_k^{3}) \leq 1 + \frac{\ell_k^{3}}{\ell_{k}^{4}}  \leq 2,\] since $\ell_k \geq 1$; combining yields the desired bound
    \[\sup_{\beta\in\T}\abs{\frac{1}{N}\sum_{n=1}^{N}w(n)e(n\beta)} \lesssim \frac{1}{p_k^{1/2}}\] as claimed.
\end{proof}

We now establish decay in every Gowers norm.

\begin{lemma}
\label{l:GowersCancellation}
    For every $k\geq 1$ and $N\in I^{(k)}$, the weight defined by \eqref{e:defiweight} satisfies \[\norm{w}_{U^s[N]} \lesssim_{s} \frac{1}{p_k^{1/2^s}}\] simultaneously for every $s\geq 2$. 
\end{lemma}

\begin{proof}
    Once again, let $1\leq J \leq 2\ell_k$ so that $N\in I^{(k)}_{J}$ and split \[\norm{w}_{U^s[N]} \leq \norm{w \cdot \one_{[B_{k-1}]}}_{U^s[N]} + \sum_{j=1}^{J-1}\norm{w\cdot \one_{I^{(k)}_{j}}}_{U^s[N]} + \norm{w\cdot \one_{I^{(k)}_{J}}}_{U^s[N]}.\] We look at each contribution separately. First, since $w$ is $1$-bounded, we can control \[\norm{w \cdot \one_{[B_{k-1}]}}_{U^s[N]} \leq \left(\frac{B_{k-1}}{N}\right)^{\frac{s+1}{2^s}} \leq \left(\frac{1}{\ell_k}\right)^{\frac{3(s+1)}{2^s}} \leq \frac{1}{p_{k}^{1/2^s}},\] where we argue as above for the final two inequalities. For the full blocks $I_{j}^{(k)}$ we again observe that \[w(n) = \eps^{(k)}_{j}V_k(n+a^{(k)}_{j})\] on these blocks and use  \Cref{c:GowersNormsCharacterWord}. Indeed, since \[\abs{I_{j}^{(k)}} = A_k(k+1)^{j-1}k \geq A_k \geq \ell_k^{2},\] we can bound \[\norm{w\cdot \one_{I^{(k)}_{j}}}_{U^s[N]} \lesssim_s \frac{1}{p_{k}^{1/2^s}}\left(\frac{\abs{I_{j}^{(k)}}}{N}\right)^{\frac{s+1}{2^s}}.\] Now, since $N\in I_{J}^{(k)}$, in particular, $N\geq A_k(k+1)^{J-1}$, and hence \[\frac{\abs{I_{j}^{(k)}}}{N}  \leq \frac{k}{(k+1)^{J-j}} \leq \frac{1}{(k+1)^{J-j-1}}.\] It follows that \[\sum_{j=1}^{J-1}\norm{w\cdot \one_{I^{(k)}_{j}}}_{U^s[N]} \lesssim_s \frac{1}{p_{k}^{1/2^s}}\sum_{j=1}^{J-1} \left((k+1)^{\frac{s+1}{2^s}}\right)^{-(J-1-j)}\lesssim_{s} \frac{1}{p_{k}^{1/2^s}}\] after estimating the geometric sum. For the remaining contribution abbreviate \[I \coloneqq I_{J}^{(k)} \cap [N],\] so that $I \subseteq [N]$ and \[\norm{w\cdot \one_{I^{(k)}_{J}}}_{U^s[N]} = \norm{w\cdot \one_{I}}_{U^s[N]}.\] If $\abs{I} \geq \ell_k^2$, then by \Cref{c:GowersNormsCharacterWord} as above \[\norm{w \cdot \one_{I}}_{U^s[N]} \lesssim_s \frac{1}{p_k^{1/2^s}}\left(\frac{\abs{I}}{N}\right)^{\frac{s+1}{2^s}} \leq \frac{1}{p_k^{1/2^s}}\] and we are done. Otherwise, simply use that $w$ is $1$-bounded to bound \[\norm{w \cdot \one_{I}}_{U^s[N]} \leq \left(\frac{\abs{I}}{N}\right)^{\frac{s+1}{2^s}} \leq \left(\frac{\ell_k^{2}}{N}\right)^{\frac{s+1}{2^s}} \leq \left(\frac{\ell_k^{2}}{\ell_k^{4}}\right)^{\frac{s+1}{2^s}} \leq \frac{1}{p_k^{1/2^s}},\] where we use that $N\geq A_k = \ell_k^4$.
\end{proof}

Putting everything together we get the concrete rate of \Cref{t:MainResult}.

\begin{proposition}
\label{p:ConcreteRate}
    Let $w$ be defined by \eqref{e:defiweight}. Then for every $N\geq 2$ \[\sup_{\beta\in \T}\abs{\frac{1}{N}\sum_{n=1}^{N}w(n)e(n\beta)} \lesssim \left(\frac{\psi(\log N)}{\log N}\right)^{1/2},\] and  \[\norm{w}_{U^s[N]} \lesssim_{s} \left(\frac{\psi(\log N)}{\log N}\right)^{1/2^{s}},\] simultaneously for every $s\geq 2$.
\end{proposition}

\begin{proof}
    Suppose first that $N\in I^{(k)}$ for some $k\geq 1$. Then \Cref{l:FourierCancellation} and \Cref{l:GowersCancellation} yield \[\sup_{\beta\in \T}\abs{\frac{1}{N}\sum_{n=1}^{N}w(n)e(n\beta)} \lesssim \frac{1}{p_k^{1/2}},\] and \[\norm{w}_{U^s[N]} \lesssim_{s} \frac{1}{p_k^{1/2^{s}}}.\] Recall now that $\ell_k=d_kp_k$, so that \[\frac{1}{p_k} = \frac{d_k}{\ell_k},\] and by construction \[d_k \leq \frac{\psi(\log(d_kp_k))}{1+\log k},\] so we have \[\frac{1}{p_k}\leq \frac{\psi(\log(d_kp_k))}{\ell_k(1+\log k)}.\] Since $N\leq B_k = \ell_k^{4}(k+1)^{2\ell_k}$, also \[\log N \leq \log B_k = 4\log \ell_k +2\ell_{k}\log (k+1) \lesssim \ell_k(1+\log k),\] and hence \[\frac{1}{p_k}\lesssim \frac{\psi(\log(d_kp_k))}{\log N}.\] Finally, since $N\geq A_k = \ell_k^{4}$, we have \[\log N \geq \log A_k \geq \log \ell_k,\] and hence, since $\psi$ is non-decreasing and $d_kp_k = \ell_k$, \[\frac{1}{p_k}\lesssim \frac{\psi(\log(N))}{\log N},\] which directly implies the result. If now $N$ falls in one of the gaps, $B_k < N \leq A_{k+1}$, then since $w\equiv 0$ on $(B_k,A_{k+1}]\cap \Z$, for the largest Fourier coefficient we have \[\sup_{\beta\in \T}\abs{\frac{1}{N}\sum_{n=1}^{N}w(n)e(n\beta)} = \frac{B_k}{N}\sup_{\beta\in \T}\abs{\frac{1}{B_k}\sum_{n=1}^{B_k}w(n)e(n\beta)} \lesssim \frac{B_k}{N} \left(\frac{\psi(\log B_k)}{\log B_k}\right)^{1/2}.\] And since $\psi$ is non-decreasing we obtain 
    \begin{align}
    \frac{B_k}{N} \left(\frac{\psi(\log B_k)}{\log B_k}\right)^{1/2} &\leq \frac{B_k}{N} \cdot \left(\frac{\psi(\log N)}{\log B_k}\right)^{1/2} \\
    & \leq \left(\frac{\psi(\log N)}{\log N}\right)^{1/2} \cdot \left( \frac{B_k (\log N)^{1/2}}{N (\log B_k)^{1/2}} \right) \lesssim \left(\frac{\psi(\log N)}{\log N}\right)^{1/2},\end{align}
    where we use that \[\frac{B_k (\log N)^{1/2}}{N (\log B_k)^{1/2}} \lesssim 1,\] because $B_k\leq N$. This shows that \[\sup_{\beta\in \T}\abs{\frac{1}{N}\sum_{n=1}^{N}w(n)e(n\beta)} \lesssim \left(\frac{\psi(\log N)}{\log N}\right)^{1/2}.\] For the Gowers norms, we can bound \[\norm{w}_{U^s[N]} = \left(\frac{B_k}{N}\right)^{\frac{s+1}{2^s}}\norm{w}_{U^s[B_k]} \lesssim_s \left(\frac{B_k}{N}\right)^{\frac{s+1}{2^s}}\left(\frac{\psi(\log B_k)}{\log B_k}\right)^{1/2^s},\] and arguing as above along with using that \[\left(\frac{B_k}{N}\right)^{s+1}\left(\frac{\log N}{\log B_k}\right) \leq 1\] for $B_k \leq N$ yields \[\norm{w}_{U^s[N]} \lesssim_{s} \left(\frac{\psi(\log N)}{\log N}\right)^{1/2^{s}}.\]
\end{proof}

\subsection{Pointwise Failure}

We now construct a system $(X,\mu,T)$ and a function $f\in L^{\infty}(\mu)$ so that the ergodic averages weighted by \eqref{e:defiweight} fail to converge on a set of positive measure. In fact, we show the following stronger divergence property.

\begin{proposition}
\label{p:PointwiseFailure}
    Let $w$ be the weight defined by \eqref{e:defiweight}. Then, there exists a measure preserving system $(X,\mu,T)$ and a measurable function $f\colon X\to \set{-1,1}$ so that \[\liminf_{N\to\infty}\frac{1}{N}\sum_{n=1}^{N}w(n)f(T^nx) = -1, \qquad \limsup_{N\to\infty}\frac{1}{N}\sum_{n=1}^{N}w(n)f(T^nx) = 1\] pointwise almost everywhere.
\end{proposition}

To improve readability, we begin by describing the system and setting up the divergence mechanism. Let $X = \set{-1,1}^{\Z}$, endowed with the usual product topology and Borel $\sigma$-algebra, $T$ be the shift map $(x_n)_{n\in\Z} \mapsto (x_{n+1})_{n\in \Z}$ and $f\colon X\to \set{-1,1}$ be the projection onto the $0$-th coordinate $f(x)=x_0$. For each $k\geq 1$ and $(\eps,a) \in \EE_k$, see \eqref{e:signtrans}, let $x^{(k,\eps,a)} \in X$ be the point defined by \[x^{(k,\eps,a)}_n = \eps V_{k}(n+a),\ n\in\Z,\] and define a probability measure $\mu_k$ on $X$ by \[\mu_k \coloneqq \frac{1}{2\ell_{k}}\sum_{(\eps,a)\in \EE_k} \delta_{x^{(k,\eps,a)}}.\] Observe that for any $n\in \Z$, \[(Tx^{(k,\eps,a)})_{n} = x^{(k,\eps,a)}_{n+1} = \eps V_k(n+a+1) = x_{n}^{(k,\eps,a+1)},\] where we regard $a+1 \mod \ell_k$. It follows that for any measurable $A\subseteq X$ \[\delta_{x^{(k,\eps,a)}}(T^{-1}A) = \delta_{x^{(k,\eps,a+1)}}(A),\] and hence, since $(\eps,a)\mapsto (\eps,a+1)$ is a bijection of $\EE_k$ when regarding $a+1 \mod \ell_k$,
\begin{align}
    \mu_k(T^{-1}A) &= \frac{1}{2\ell_{k}}\sum_{(\eps,a)\in \EE_k} \delta_{x^{(k,\eps,a)}}(T^{-1}A)\\
    &= \frac{1}{2\ell_{k}}\sum_{(\eps,a)\in \EE_k} \delta_{x^{(k,\eps,a+1)}}(A)= \frac{1}{2\ell_{k}}\sum_{(\eps,a)\in \EE_k} \delta_{x^{(k,\eps,a)}}(A) = \mu_k(A),
\end{align}
so $\mu_k$ is an invariant measure for the shift.

Let $\mu$ be the limit of some weakly converging subsequence $(\mu_{k_m})_{m}$ of $(\mu_k)_{k}$. Then $\mu$ is again a shift invariant probability measure. Finally, for every $k\geq 1$ we consider the sets \[E_k \coloneqq \bigcup_{(\eps,a)\in \EE_k} \set{x\in X \colon x_n = \eps V_k(n+a),\ \forall 1\leq n\leq B_k},\] and note each is a finite union of cylinder sets, and hence clopen. Crucially, they are also eventually large.

\begin{lemma}
\label{l:SmallMeasure}
    For every $k\geq 1$ and $j\geq k+1$ we have \[\mu_{j}(\comp{E_k}) \leq \frac{B_k}{d_{k+1}} \leq 2^{-(k+1)}.\] In particular $\mu(\comp{E_k}) \leq 2^{-(k+1)}$.
\end{lemma}

\begin{proof}
    Let $j\geq k+1$. First, by \Cref{l:Restriction}, for every $m\in \Z$ there exists a sign $\eta = \eta_{j,m}\in \set{-1,1}$, such that \[V_{j}(n) = \eta V_{k}(n),\] for every $md_{k+1} \leq n < (m+1)d_{k+1}$. Observe now that for any $(\sigma,b)\in \EE_{j}$, if there exists some $m\in \Z$ so that \[\set{b+1,\ldots,b+B_k} \subseteq\set{md_{k+1},\ldots,(m+1)d_{k+1}-1},\] then for every $n\in \set{1,\ldots,B_k}$ \[x^{(j,\sigma,b)}_n = \sigma V_{j}(n+b) = \sigma\eta V_k(n+b),\] so letting $\eps=\sigma\eta\in\set{-1,1}$ and $a\in \set{0,\ldots,\ell_k-1}$ so that $a\equiv b \mod \ell_k$, we have \[x^{(j,\sigma,b)}_n = \eps V_k(n+a),\] which implies that $x^{(j,\sigma,b)}\in E_k$, and hence \[\delta_{x^{(j,\sigma,b)}}(\comp{E_k})=0.\] Taking the contrapositive of this observation, it follows that if for some $(\sigma,b)\in \EE_j$ \[\delta_{x^{(j,\sigma,b)}}(\comp{E_k})=1\] then necessarily 
    \begin{equation}
    \label{e:notcontained}
        \set{b+1,\ldots,b+B_k} \not\subseteq\set{md_{k+1},\ldots,(m+1)d_{k+1}-1},
    \end{equation}
    for every $m\in \Z$, and thus \[\mu_{j}(\comp{E_{k}}) \leq \frac{1}{2\ell_j}\abs{\set{(\sigma, b) \in \EE_j \colon \eqref{e:notcontained} \text{ holds for every }m\in \Z}}.\] We now estimate this quantity. For $b\in \set{0,\ldots,\ell_j-1}$, write \[b+1= md_{k+1}+r\] with $m\in \Z$ and $0\leq r < d_{k+1}$. It suffices to check \eqref{e:notcontained} only for this $m$; for every other value it holds immediately. Now, for this $m$, \eqref{e:notcontained} holds if and only if $r\geq d_{k+1}-B_{k}+1$. Since also $r<d_{k+1}$, this means that \[d_{k+1}-B_{k}+1 \leq r \leq d_{k+1}-1,\] so there are at most $B_k$ residues $r$, such that \eqref{e:notcontained} is satisfied by $b+1 = md_{k+1}+r$. Any one of these residues defines at most $\ell_j/d_{k+1}$ many possible $b\in \set{0,\ldots,\ell_j-1}$ via $b+1 = md_{k+1}+r$, and so there are at most $B_{k}\ell_{j}d_{k+1}^{-1}$ many $b\in \set{0,\ldots,\ell_j-1}$ so that \eqref{e:notcontained} holds for every $m\in \Z$. Since every such $b$ can appear accompanied by two possible signs, we conclude that \[\abs{\set{(\sigma, b) \in \EE_j \colon \eqref{e:notcontained} \text{ holds for every }m\in \Z}} \leq \frac{2B_{k}\ell_j}{d_{k+1}},\] and in turn \[\mu_{j}(\comp{E_{k}}) \leq \frac{2B_{k}\ell_j}{2\ell_jd_{k+1}} = \frac{B_k}{d_{k+1}} \leq 2^{-(k+1)},\] where the last inequality is just by construction. Since $\comp{E_k}$ is clopen, we conclude that \[\mu(\comp{E_k}) = \lim_{m\to \infty}\mu_{k_m}(\comp{E_k})\leq 2^{-(k+1)}.\]
\end{proof}

We are now ready to prove \Cref{p:PointwiseFailure}.

\begin{proof}[Proof of \Cref{p:PointwiseFailure}]
    Let $(X,\mu,T)$ be the measure preserving system described above and consider the set \[\Omega \coloneqq \set{x\in X \colon \text{ there exists }k_0\geq 1 \text{ so that } x\in E_k \text{ for every } k\geq k_0} = \bigcup_{k_0\geq 1} \bigcap_{k \geq k_0} E_k.\] Then, by \Cref{l:SmallMeasure} and Borel-Cantelli, $\mu(\Omega)=1$. We now show that \[\liminf_{N\to\infty}\frac{1}{N}\sum_{n=1}^{N}w(n)f(T^nx) = -1,\qquad \limsup_{N\to\infty}\frac{1}{N}\sum_{n=1}^{N}w(n)f(T^nx) = 1\] for any $x\in \Omega$, which establishes the result. We only show the statement for the limit superior, the argument for the limit inferior is analogous.
    
    Given $x\in \Omega$, we have that $x\in E_k$ for every large enough $k\geq 1$ and so there exists some $(\eps,a)\in \EE_k$, such that \[x_n = \eps V_k(n+a)\] for every $1\leq n\leq B_k$. Letting $1\leq j\leq 2\ell_k$, such that $(\eps,a) = (\eps^{(k)}_j,a^{(k)}_{j})$, for every $n\in I^{(k)}_{j}\subseteq [B_k]$, we have \[w(n)f(T^nx) = w(n)x_n = (\eps V_k(n+a))^{2} = 1.\] It follows that for $N=A_k(k+1)^{j}$ \[\frac{1}{N}\sum_{n=1}^{N}w(n)f(T^nx) = \frac{1}{N}\sum_{n=1}^{A_k(k+1)^{j-1}}w(n)f(T^nx) + \frac{\abs{I_j^{(k)}}}{N},\] and since \[\abs{\frac{1}{N}\sum_{n=1}^{A_k(k+1)^{j-1}}w(n)f(T^nx)} \leq \frac{A_k(k+1)^{j-1}}{A_k(k+1)^j} = \frac{1}{k+1},\] and \[\frac{\abs{I_j^{(k)}}}{N} = \frac{A_k(k+1)^{j-1}k}{A_k(k+1)^j} = \frac{k}{k+1},\] we get \[\frac{1}{N}\sum_{n=1}^{N}w(n)f(T^nx) \geq \frac{k}{k+1} - \frac{1}{k+1} = 1 - \frac{2}{k+1}.\] Since this holds for every large enough $k\geq 1$, we deduce that \[\limsup_{N\to \infty}\frac{1}{N}\sum_{n=1}^{N}w(n)f(T^nx) \geq 1.\] 
\end{proof}

\begin{remark}
    It is not too hard to show that the system $(X,\mu,T)$ described above is rigid along the sequence of word-lengths $(\ell_{k})_{k}$, meaning that \[\norm{g-T^{\ell_{k}}g}_{L^2(\mu)}\to 0\] for every $g\in L^2(\mu)$. In particular, this means that $\mu$ has zero measure-theoretical entropy, so the divergence above happens in a deterministic system. Furthermore, taking $Y=\supp \mu$, one can show that $(Y,T)$ is uniquely ergodic with $\mu$ being the unique shift invariant measure, so the fact that $\mu$ has zero measure-theoretical entropy in particular implies that the topological entropy of $(Y,T)$ is zero. We leave the details to the interested reader.
\end{remark}

\begin{proof}[Proof of \Cref{t:MainResult}]
    For any fixed $\psi \colon \R_{+}\to \R_{+}$ as in the theorem, the construction above along with \Cref{p:ConcreteRate} yield the desired decay rate, and \Cref{p:PointwiseFailure} shows the pointwise divergence.
\end{proof}

\begin{remark}
    In line with \Cref{remark:summability}, one can obtain a slightly more general construction and prove the existence of weights $w\colon \Z\to \set{-1,0,1}$ with arbitrarily slowly diverging sums that nevertheless fail to be pointwise good. Concretely, given any non-decreasing $\Psi\colon \R_{+}\to \R_{+}$ with $\Psi(x) \to \infty$ as $x\to \infty$, by choosing the primes $(p_k)_{k}$ so that \[\Psi(A_k) \geq \sum_{j=1}^{k}(1+d_j\log(j+1)),\] one can obtain a weight so that \[\sum_{N=1}^{M}\frac{\norm{w}_{U^s[N]}^{2^s}}{N} \lesssim_s \Psi(M),\qquad \sum_{N=1}^{\infty}\frac{\norm{w}_{U^s[N]}^{2^s}}{N} = \infty\] simultaneously for every $s\geq 2$, for which the same almost-everywhere divergence conclusion as in \Cref{t:MainResult} holds.
\end{remark}

\section{Genericity}
\label{sec:Genericity}

In this section we show that the behavior of the weight constructed in the previous section is not isolated and is typical in a very concrete sense. The results in this section are stated and proved for weights taking values in $\set{-1,0,1}$ or the unit disk $\overline{\D}$, but hold in full generality for any bounded subset of the complex numbers containing a non-zero scalar copy of $\set{-1,0,1}$. We begin by setting the stage.

\subsection{Failure of Pointwise Goodness}

Given a bounded set $K\subseteq \C$ we let $\WW(K)$ denote the space of weights \[w\colon \Z \to K\] and recall the Besicovitch distance \[d_B(u,v)\coloneqq \limsup_{N\to\infty}\frac{1}{N}\sum_{n=1}^{N}\abs{u(n)-v(n)}.\] This is the natural notion of distance tied to pointwise convergence of bounded functions as witnessed by the following elementary inequality
\begin{equation}
\label{e:AveragesBesicovitch}
    \limsup_{N\to\infty}\norm{\frac{1}{N}\sum_{n=1}^{N}u(n)f(T^nx) - \frac{1}{N}\sum_{n=1}^{N}v(n)f(T^nx)}_{L^{\infty}(\mu)} \leq \norm{f}_{L^{\infty}(\mu)}d_B(u,v).
\end{equation}
It is, however, only a pseudo-metric; it is non-negative, symmetric, satisfies the triangle inequality and $d_B(u,u)=0$, but $d_B(u,v)=0$ need not imply $u=v$. To turn it into a metric, we identify $u\sim v$ whenever $d_B(u,v)=0$ and let $W_B(K)$ be the resulting space. It is standard that $(W_B(K),d_B)$ is a complete metric space (see e.g. \cite{LesigneQuasRosenblattWierdl}). Furthermore, \eqref{e:AveragesBesicovitch} implies that pointwise goodness is a well-defined concept for equivalence classes $w\in W_B(K)$. We begin by showing that failure of pointwise goodness is generic in this topology.

\begin{theorem}
\label{t:WeightGeneric}
    Let $K=\set{-1,0,1}$ or $K=\overline{\D}$. Then, the subset of $W_B(K)$ of weights that are not pointwise good is open and dense for the Besicovitch metric.
\end{theorem}

\Cref{t:WeightGeneric} follows quickly from the following two lemmas, whose proofs we defer to after the main argument.

\begin{lemma}
\label{l:Openness}
    Let $K=\overline{\D}$ or $K=\set{-1,0,1}$ and let $w\in W_{B}(K)$ be a weight that is not pointwise good. Then, there exists some $\eps>0$ such that every $v\in W_B(K)$ with $d_B(w,v) < \eps$ is also not pointwise good.
\end{lemma}

\begin{lemma}
\label{l:Suspension}
    Let $K=\set{-1,0,1}$ or $K=\overline{\D}$ and let $w\in W_{B}(K)$ be a weight that is not pointwise good. Then, for any $v\in W_B(K)$ and every integer $q\geq 2$ the weight $u = u_q\in W_{B}(K)$ defined by \[u(n) = \begin{cases}
        w(n/q),\ &\text{ if }q\mid n,\\
        v(n),\ &\text{ otherwise,}
    \end{cases}\]
    is not pointwise good and $d_{B}(v,u) \leq 2/q$.
\end{lemma}

\begin{proof}[Proof of \Cref{t:WeightGeneric} assuming \Cref{l:Openness,l:Suspension}]
    The fact that the subset of the space $W_{B}(K)$ of weights which are not pointwise good is non-empty is elementary; it suffices, for example, to take the indicator of a subset of the integers that doesn't admit a density. The openness of this set then follows immediately from \Cref{l:Openness} and its density is a consequence of \Cref{l:Suspension} after taking $q$ large enough to have a $d_{B}$-ball of radius $2/q$ fit inside an open set $V\subseteq W_{B}(K)$.
\end{proof}

We now prove our lemmas.
\begin{proof}[Proof of \Cref{l:Openness}]
    Since $w$ is not pointwise good, there exists some measure-preserving system $(X,\mu,T)$ and some $f\in L^{\infty}(\mu)$ such that the averages $\frac{1}{N}\sum_{n=1}^{N}w(n)f(T^nx)$ fail to converge on a set of positive measure. In particular, the averages fail to form a Cauchy sequence, so there is a set of positive measure, such that for every $x$ in this set \[\limsup_{N,M \to \infty}\abs{\frac{1}{N}\sum_{n=1}^{N}w(n)f(T^nx)-\frac{1}{M}\sum_{n=1}^{M}w(n)f(T^nx)} > 0,\] and furthermore, by continuity of the measure, we can find some $\delta>0$ such that 
    \begin{equation}
    \label{e:opennes01}
        \limsup_{N,M \to \infty}\abs{\frac{1}{N}\sum_{n=1}^{N}w(n)f(T^nx)-\frac{1}{M}\sum_{n=1}^{M}w(n)f(T^nx)} > \delta
    \end{equation}
    on a set of positive measure. For any $x$ in this set by the triangle inequality, \eqref{e:opennes01} and \eqref{e:AveragesBesicovitch} we now have \[\limsup_{N,M\to \infty} \abs{\frac{1}{N}\sum_{n=1}^{N}v(n)f(T^nx) - \frac{1}{M}\sum_{n=1}^{M}v(n)f(T^nx)} \geq \delta - 2\norm{f}_{L^{\infty}(\mu)}\cdot d_{B}(w,v),\] so choosing \[\eps \leq \frac{\delta}{2\norm{f}_{L^{\infty}(\mu)}},\] we have \[\limsup_{N,M\to \infty} \abs{\frac{1}{N}\sum_{n=1}^{N}v(n)f(T^nx) - \frac{1}{M}\sum_{n=1}^{M}v(n)f(T^nx)} > 0\] for any $v$ with $d_B(w,v)< \eps$. By the above, any such $v$ is not pointwise good.
\end{proof}
\begin{proof}[Proof of \Cref{l:Suspension}]
    Begin by observing that $u\in W_{B}(K)$ since it takes values exactly in the set where $w$ and $v$ do. Now, since $w$ is not pointwise good, there exists some measure-preserving system $(X,\mu,T)$ and some $f\in L^{\infty}(\mu)$ such that the ergodic averages of $f$ weighted by $w$ fail to converge on a set of positive measure. Consider the finite rotation $r\mapsto r+1$ on $\Z_{q}$ endowed with the Haar measure $\lambda_q$ and define the measure preserving system $(Y,\nu,S)$ as the suspension of $T$ of constant height $q$; this is, \[Y = X \times \Z_{q}, \qquad \nu = \mu \otimes \lambda_q, \qquad S(x,r) \coloneqq \begin{cases}
        (x,r+1),\ &\text{ if } r\neq q-1,\\
        (Tx,0),\ &\text{ if }r=q-1.
    \end{cases}\]
    Putting $F(x,r)=f(x)\one_{\set{0}}(r)$, we have that $F\in L^{\infty}(\nu)$, and moreover \[F(S^n(x,0)) = 0, \text{ if }q\nmid n,\] and \[F(S^{qm}(x,0)) = f(T^mx),\] otherwise. It follows that \[\frac{1}{qN}\sum_{n=1}^{qN}u(n)F(S^{n}(x,0)) = \frac{1}{q}\cdot \frac{1}{N}\sum_{m=1}^{N}w(m)f(T^mx).\] Since the right hand side above fails to converge on some set $E\subseteq X$ with $\mu(E)>0$, so does the left hand side on the set $E\times \set{0}$ of measure \[\nu(E\times \set{0}) = \frac{\mu(E)}{q} > 0.\] So $u$ is not pointwise good. The second statement follows from the fact that $u$ and $v$ coincide everywhere but on the arithmetic progression $q\Z$, where we can simply bound $\abs{v(n)-u(n)}\leq 2$.
\end{proof}

\subsection{Failure of Pointwise Goodness under Fourier and Gowers Decay}

We now show that the property of not being pointwise good remains typical even when restricting to weights with decaying Gowers norms. Concretely, if we consider the subspace $U_B(K)\subseteq W_B(K)$ defined by \[U_B(K) \coloneqq \set{w\in W_{B}(K) \colon \norm{w}_{U^s[N]} \to 0,\ \forall s\geq 2},\] we have the following.

\begin{cor}
\label{c:UniformWeightGeneric}
    Let $K=\set{-1,0,1}$ or $K=\overline{\D}$. Then, the subset of $U_B(K)$ of weights that are not pointwise good is open and dense for the Besicovitch metric.
\end{cor}

We first justify that the space $U_B(K)$ is well-defined and that it is in fact a closed subspace of $W_B(K)$. It all boils down to the following inequality, which we state in the broad setting of bounded $K$.

\begin{lemma}
\label{l:Inequalities}
    Let $K\subseteq \C$ be bounded and let $u,v\in \WW(K)$. Then, for every $s\geq 2$ \[\limsup_{N\to \infty}\abs{\norm{u}_{U^s[N]} - \norm{v}_{U^s[N]}} \lesssim_{s,K} d_B(u,v)^{1/2^s}.\] 
\end{lemma}

\begin{proof}
    By the reverse triangle inequality we can bound \[\abs{\norm{u}_{U^s[N]} - \norm{v}_{U^s[N]}} \leq \norm{u-v}_{U^s[N]},\] so after expanding \[\norm{u-v}_{U^s[N]}^{2^s} = \frac{1}{N^{s+1}}\sum_{n,h_1,\ldots,h_s \in \Z} \prod_{\omega\in \set{0,1}^{s}} \CC^{\abs{\omega}}(u(n+\omega \cdot h)- v(n+\omega\cdot h))\one_{[N]}(n+\omega\cdot h),\] taking absolute values, factoring out the $\omega = \mathbf{0} := (0,\dots,0) \in \{0,1\}^s$ term and bounding all the remaining factors by $\mathrm{diam}(K)\lesssim_{K}1$ one gets \[\norm{u-v}_{U^s[N]}^{2^s} \lesssim_{s,K} \frac{1}{N}\sum_{n=1}^{N}\abs{u(n)-v(n)}\cdot \frac{1}{N^s}\sum_{h_1,\ldots,h_s\in \Z}\prod_{\omega \in \set{0,1}^{s}\setminus\set{\0}}\one_{[N]}(n+\omega \cdot h).\] One can then simply bound \[\frac{1}{N^s}\sum_{h_1,\ldots,h_s\in \Z}\prod_{\omega \in \set{0,1}^{s}\setminus\set{\0}}\one_{[N]}(n+\omega \cdot h) \leq 2^{2^s-1}\frac{1}{N^s}\sum_{\abs{h_1},\ldots,\abs{h_s}\leq N} \lesssim_s 1,\] and the result follows.
\end{proof}

As a consequence we get that the expression \[\norm{w}_{U^s[N]} \to 0\] is well defined when regarding $w$ as an equivalence class in $W_B(K)$. In turn, this implies that the subspace $U_B(K)$ is well defined. Moreover, again using \Cref{l:Inequalities} it is not hard to see that it is in fact closed, and so $(U_B(K),d_B)$ is again a complete metric space. We now prove the corollary.

\begin{proof}[Proof of \Cref{c:UniformWeightGeneric}]
    Let $\BB$ be the subset of $W_B(K)$ of weights that are not pointwise good. Then the corresponding subset of $U_B(K)$ is just $U_B(K)\cap \BB$. By \Cref{t:WeightGeneric}, $\BB$ is open in $W_B(K)$, and hence so is $U_B(K)\cap \BB$ in $U_B(K)$; by \Cref{t:MainResult} it is crucially non-empty. It remains to show that it is dense, and to do this it suffices to show that the weights constructed in \Cref{l:Suspension} remain in $U_B(K)$ if $w,v\in U_B(K)$. Indeed, observe that given any $q\geq 2$, the weight $u=u_q$ furnished by \Cref{l:Suspension} can be written as \[u(n) = v(n) - \one_{q\mid n}v(n) + \one_{q\mid n}w(n/q),\] and since by orthogonality \[\frac{1}{q}\sum_{r=0}^{q-1}e\left(\frac{rn}{q}\right) = \one_{q\mid n},\] we can actually express \[u(n) = v(n) - \frac{1}{q}\sum_{r=0}^{q-1}v(n)e\left(\frac{rn}{q}\right) + \one_{q\mid n}w(n/q).\] By the triangle inequality for Gowers norms, we then have \[\norm{u}_{U^s[N]} \leq \norm{v}_{U^s[N]} + \frac{1}{q}\sum_{r=0}^{q-1}\norm{ve(\cdot r/q)}_{U^s[N]} + \norm{w(\cdot/q) \one_{q\mid \cdot}}_{U^s[N]},\] and by their modulation invariance for $s\geq 2$, \[\frac{1}{q}\sum_{r=0}^{q-1}\norm{ve(\cdot r/q)}_{U^s[N]} = \frac{1}{q}\sum_{r=0}^{q-1}\norm{v}_{U^s[N]} = \norm{v}_{U^s[N]}.\] If we now expand \[\norm{w(\cdot/q) \one_{q\mid \cdot}}_{U^s[N]}^{2^s} = \frac{1}{N^{s+1}}\sum_{n,h_1,\ldots,h_s\in \Z}\prod_{\omega\in \set{0,1}^{s}}\CC^{\abs{\omega}}w\left(\frac{n+\omega\cdot h}{q}\right)\one_{q\mid n+\omega\cdot h}\one_{[N]}(n+\omega\cdot h),\] we can observe that for a summand to be non-zero we need $q\mid n+\omega \cdot h$ simultaneously for every $\omega\in\set{0,1}^{s}$, so necessarily $n,h_1,\ldots,h_s\in q\Z$. Writing $n=qm$ and $h_i=q\ell_i$, we can now express \[\norm{w(\cdot/q) \one_{q\mid \cdot}}_{U^s[N]}^{2^s} = \frac{1}{N^{s+1}}\sum_{m,\ell_1,\ldots,\ell_s\in \Z}\prod_{\omega\in \set{0,1}^{s}}\CC^{\abs{\omega}}w(m+\omega \cdot \ell)\one_{[N]}(q(m+\omega\cdot \ell)),\] and since $\one_{[N]}(q(m+\omega\cdot \ell)) = \one_{[\lfloor N/q \rfloor]}(m+\omega\cdot \ell)$, we have
    \begin{align}
        \norm{w(\cdot/q) \one_{q\mid \cdot}}_{U^s[N]}^{2^s} &= \frac{1}{N^{s+1}}\sum_{m,\ell_1,\ldots,\ell_s\in \Z}\prod_{\omega\in \set{0,1}^{s}}\CC^{\abs{\omega}}w(m+\omega \cdot \ell)\one_{[\lfloor N/q \rfloor]}(m+\omega\cdot \ell)\\
        &= \left(\frac{\lfloor N/q \rfloor}{N}\right)^{s+1} \norm{w}_{U^s[\lfloor N/q \rfloor]}^{2^s},
    \end{align}
    where we assume that $N\geq q$. This means that \[\norm{u}_{U^s[N]} \leq 2\norm{v}_{U^s[N]} + \left(\frac{\lfloor N/q \rfloor}{N}\right)^{\frac{s+1}{2^s}} \norm{w}_{U^s[\lfloor N/q \rfloor]},\] so since \[\norm{w}_{U^s[N]}, \norm{v}_{U^s[N]}\to 0, \qquad \left(\frac{\lfloor N/q \rfloor}{N}\right)^{\frac{s+1}{2^s}} \to \left(\frac{1}{q}\right)^{\frac{s+1}{2^s}},\] we get that \[\norm{u}_{U^s[N]}\to 0.\]
\end{proof}

In fact, the argument above also shows that failure of pointwise goodness remains generic even when we impose quantitative decay arbitrarily close to the thresholds of \Cref{t:PositiveResultFourier} and \Cref{p:PositiveResult} that ensure pointwise goodness. We must be a bit careful with how we phrase this, since two members of an equivalence class in $W_{B}(K)$ need not have the same finite-scale behavior. To overcome this, we state this in terms of equivalence classes admitting a representative with the prescribed decay rates. For simplicity denote \[\Phi_{N}(w)\coloneqq \sup_{\beta\in \T}\abs{\frac{1}{N}\sum_{n=1}^{N}w(n)e(n\beta)}\] and given a non-decreasing function $\psi\colon \R_{+}\to \R_{+}$ with $\psi(x)\to \infty$ as $x\to \infty$, let $U_{B}^{\psi}(K)$ be the set of equivalence classes in $W_{B}(K)$ admitting a representative $w\in \WW(K)$ such that for every $N\geq 2$ \[\Phi_{N}(w) \lesssim \left(\frac{\psi(\log N)}{\log N}\right)^{1/2},\qquad \norm{w}_{U^s[N]} \lesssim_{s} \left(\frac{\psi(\log N)}{\log N}\right)^{1/2^s},\ \forall s\geq 2.\]

\begin{cor}
\label{c:PrescribedDecayGeneric}
    Let $K=\set{-1,0,1}$ or $K=\overline{\D}$. Then, for any non-decreasing $\psi\colon \R_{+}\to \R_{+}$ with $\psi(x)\to \infty$ the subset of $U_{B}^{\psi}(K)$ of weights that are not pointwise good is open and dense for the Besicovitch metric.
\end{cor}

\begin{proof}
    We just need to check that for each $q\geq 2$ the approximant $u=u_q$ constructed from $w,v \in U_{B}^{\psi}(K)$, which we regard as weights with the prescribed decay rather than equivalence classes, also remains in $U_{B}^{\psi}(K)$. Looking at the proof of \Cref{c:UniformWeightGeneric}, we have \[\norm{u}_{U^s[N]} \lesssim \norm{v}_{U^s[N]} +  \norm{w}_{U^s[\lfloor N/q \rfloor]},\] and arguing similarly \[\Phi_N(u) \lesssim \Phi_N(v) + \Phi_{\lfloor N/q \rfloor}(w).\] The first summands in both inequalities readily decay as prescribed. For the second ones, note that for any fixed $q\geq 2$ and $N\geq 2q$ we have \[\frac{1}{\log \lfloor N/q \rfloor} \lesssim_{q} \frac{1}{\log N},\] and thus, since $\psi$ is non-decreasing \[\frac{\psi(\log \lfloor N/q \rfloor)}{\log \lfloor N/q \rfloor} \lesssim_{q} \frac{\psi(\log N)}{\log N},\] which shows that \[\Phi_{\lfloor N/q \rfloor}(w) \lesssim_{q} \left(\frac{\psi(\log N)}{\log N}\right)^{1/2}, \qquad \norm{w}_{U^s[\lfloor N/q \rfloor]} \lesssim_{s,q} \left(\frac{\psi(\log N)}{\log N}\right)^{1/2^s},\ \forall s\geq 2\] as desired.
\end{proof}

\section{Probabilistic Balayage and Gowers Uniform Sets}
\label{sec:GowersUniform}

This section is devoted to proving \Cref{t:GowersUniformBad}. We begin by putting this result into context. For $s\geq 1$, we say that a set $E\subseteq \Z$ is \textit{Gowers uniform of order $s$} if $d(E)>0$ and \[\lim_{N\to \infty}\norm{\one_E - d(E)}_{U^s[N]} = 0,\] and that it is \textit{Gowers uniform of all orders} or simply \textit{Gowers uniform} if it is Gowers uniform of order $s$ for every $s\geq 1$. 

Examples of Gowers uniform sets include level sets of the prime $\Omega$ function modulo $q\geq 2$ \[\set{n\in \N \colon \Omega(n) \equiv a \mod q}\] and arithmetic Bohr-type sets of the form \[\set{n\in \N \colon \alpha \Omega(n) \mod 1 \in U}\] for $\alpha\notin \Q$ and $U\subseteq \T$ with $m_{\T}(U)>0$ and $m_{\T}(\partial U) = 0$, where $m_{\T}$ denotes the Haar measure on $\T$; see \cite{FrantzikinakisHost2017}. In the context of \textit{automatic sequences}, J. Konieczny \cite{Konieczny} showed that the set of $1$s of the Thue-Morse sequence is Gowers uniform, and more generally, Fan and Konieczny \cite{FanKonieczny} showed that for any base $b$ and any $q\geq 2$ coprime with $b-1$, so is the set \[\set{n\in \N \colon s_{b}(n) \equiv a \mod q}\] where $s_{b}(n)$ sums the digits of $n$ in base $b$. Moreover, it is not hard to check that a random set obtained by choosing each integer independently with probability $\delta \in (0,1)$ will almost surely be Gowers uniform.

These sets enjoy a great deal of combinatorial and dynamical properties. In the infinitary setting, Frantzikinakis and Host \cite{FrantzikinakisHost2017} showed that Gowers uniform sets are not only norm good but even sets of polynomial multiple mean convergence and polynomial multiple recurrence, with a refinement by A. Koutsogiannis \cite{Koutsogiannis} extending recurrence to closest integer polynomial iterates and mean convergence to integer parts of real polynomials. Meanwhile, in the finitary setting, Gowers uniform sets were shown to have the expected counts of affine-linear configurations of finite complexity in the first $N$ integers by B. Green and T. Tao \cite{GreenTao}, with polynomial analogues due to Tao and T. Ziegler \cite{TaoZiegler}, S. Peluse \cite{Peluse}, and L. Matthiesen, Teräväinen and M. Wang \cite{MatthiesenTeravainenWang} under suitable conditions. More recently, T. Radi\'c \cite{Radic} and F. Hern\'andez and Radi\'c \cite{HernandezRadic} showed that these sets can additionally have rich infinite sumset structure if one replaces Gowers norms in the definition above for the stronger \textit{local uniformity seminorms} introduced by Host and Kra \cite{HostKra}.

In contrast, \Cref{t:GowersUniformBad} says that such sets can nevertheless fail to be pointwise good, regardless of how dense they are.\\

The proof of \Cref{t:GowersUniformBad} will make use of the probabilistic balayage technique developed by Bourgain \cite{BourgainUnpublished} to construct a random set $E$ by independently choosing integers with probability given by the weight $w$ we've constructed. Here is the framework:

Fixing $\eps\in (0,1/2]$ and letting $a\colon \N\to \R$ be a $1$-bounded sequence, put \[\sigma_n = (1-\eps) + \eps a(n)\in [0,1],\] and construct a random set $E\subseteq \N$ by independently setting \[\P(n\in E) = \sigma_n.\] Alternatively, let $(X_n)_{n\geq 1}$ be independent Bernoulli random variables with parameter $\sigma_n$ and put \[E = \set{n\in \N \colon X_n = 1}.\] Under this perspective, we show the following proposition, which directly implies \Cref{t:GowersUniformBad}.

\begin{proposition}
\label{p:balayage}
    Let $a\colon \N\to \R$ be a $1$-bounded sequence that is not pointwise good for $L^{\infty}$ functions and satisfies $\norm{a}_{U^s[N]}\to 0$ for some $s\geq 2$. Then, almost surely, a set $E$ constructed by probabilistic balayage has $d(E)=1-\eps$, is Gowers uniform of order $s$, and is not pointwise good for $L^{\infty}$ functions.
\end{proposition}

\begin{proof}[Proof of \Cref{t:GowersUniformBad} assuming \Cref{p:balayage}]
    Let $w$ be a weight given by \Cref{t:MainResult} for any function $\psi$ so that $\norm{w}_{U^s[N]}\to 0$ for every $s\geq 2$. Applying \Cref{p:balayage} to the restriction to $\N$ of $w$ and intersecting the countably many full probability events obtained (one for each $s \geq 2$) yields the result.
\end{proof}

To improve readability, we split the proof of \Cref{p:balayage} into several steps. To begin with, define the centered random variables $Y_n = X_n - \sigma_n$ and observe that $\abs{Y_n}\leq 1$, $\E(Y_n) = 0$ and $\E(\abs{Y_n}^2) = \sigma_n(1-\sigma_n) \leq \sigma_n$. This allows us to write 
\begin{equation}
\label{e:usefulexpressionYn}
    \one_E(n) - (1-\eps) = Y_n + \eps a(n).
\end{equation}
As the next result shows, Gowers uniformity of any order is enough to ensure that almost surely $d(E)=1-\eps$.

\begin{lemma}
\label{l:densityan}
    Suppose $\norm{a}_{U^s[N]}\to 0$ for some $s\geq 1$. Then $d(E)=1-\eps$ almost surely.
\end{lemma}

\begin{proof}
    Simply write \[\frac{\abs{E \cap [N]}}{N} = \frac{1}{N}\sum_{n=1}^{N}\one_{E}(n) = \frac{1}{N}\sum_{n=1}^{N}Y_n + \eps\cdot\frac{1}{N}\sum_{n=1}^{N}a(n) + (1-\eps).\] The first average converges to $0$ as $N\to \infty$ by Kolmogorov's Strong Law and so does the second one since \[\abs{\frac{1}{N}\sum_{n=1}^{N}a(n)} = \norm{a}_{U^1[N]} \lesssim_{s} \norm{a}_{U^s[N]}\to 0,\] thus only leaving the constant term $(1-\eps)$.
\end{proof}

We now show that the random perturbation $Y_n$ doesn't break uniformity. 

\begin{lemma}
\label{l:Ynuniform}
    Almost surely, for every $s\geq 1$ one has $\norm{Y}_{U^s[N]} \lesssim_{s} N^{-1/2^{s+2}}$ for every sufficiently large $N\geq 1$.
\end{lemma}

\begin{proof}
    It suffices to show that for every fixed $s\geq 1$, one has $\norm{Y}_{U^s[N]} \lesssim N^{-1/2^{s+2}}$ almost surely, and take a countable intersection. Start by observing that \[\E\left(\norm{Y}_{U^s[N]}^{2^s}\right) = \frac{1}{N^{s+1}}\sum_{\substack{n,h_1,\ldots,h_s \in \Z\\ n+\omega \cdot h \in [N],\ \forall \omega \in \set{0,1}^s}}\E\left(\prod_{\omega\in\set{0,1}^s}Y_{n+\omega\cdot h}\right).\] Given $n,h_1,\ldots,h_s$, if all $n+\omega \cdot h$ are distinct, then, by independence and the fact that each $Y_n$ has zero expectation, the corresponding summand is $0$. It follows that the only surviving contributions to the sum come from $n,h_1,\ldots,h_s$, such that $n+\omega \cdot h = n+\omega'\cdot h$ for some $\omega'\neq \omega$, in which case one coordinate is determined from all of the others. Since $n\in [N]$ and $\abs{h_i}\leq N$, it follows that there are $\leq C_s N^s$ such cases, for some constant $C_s>0$, so bounding their expectation by $1$ we get 
    \begin{equation}
    \label{e:expectationGowersYn}
        \E\left(\norm{Y}_{U^s[N]}^{2^s}\right) \leq \frac{C_sN^s}{N^{s+1}} \leq \frac{C_s}{N}.
    \end{equation}
    By Markov's Inequality along the sequence of squares, we get\[\P\left(\norm{Y}_{U^s[N^2]}^{2^s} > \frac{C_s}{N^{1/2}}\right) \leq \frac{1}{N^{3/2}},\] which is summable in $N$, so by Borel-Cantelli, we get that \[\norm{Y}^{2^{s}}_{U^s[N^2]} \lesssim_{s} N^{-1/2} = (N^2)^{-1/4}\] almost surely. To pass to the full sequence, we appeal to an appropriate H\"{o}lder continuity of Gowers norms for $1$-bounded functions. Specifically, given $N\geq 1$, let $M=\lfloor N^{1/2}\rfloor$ so that $M^2 \leq N < (M+1)^{2}$ and split \[\norm{Y}_{U^s[N]} \leq \norm{Y\one_{[M^2]}}_{U^s[N]} + \norm{Y\one_{(M^2,N]}}_{U^s[N]}.\] For the first term, by the above we have
    \begin{align}
        \norm{Y\one_{[M^2]}}_{U^s[N]} = \left(\frac{M^2}{N}\right)^{\frac{s+1}{2^s}} \norm{Y}_{U^s[M^2]} &\lesssim_{s} \left(\frac{M^2}{N}\right)^{\frac{s+1}{2^s}}  M^{-\frac{1}{2^{s+1}}}\\
        &=\left(\frac{M^2}{N}\right)^{\frac{s+1}{2^s} - \frac{1}{2^{s+2}}}N^{-1/2^{s+2}} \leq N^{-1/2^{s+2}},
    \end{align}
    where we use that $M^2/N \leq 1$ and that $(s+1)/2^s - 1/2^{s+2} > 0$. Since $Y$ is $1$-bounded, for the second term we can simply bound \[\norm{Y\one_{(M^2,N]}}_{U^s[N]} \leq \left(\frac{N-M^{2}}{N}\right)^{\frac{s+1}{2^s}},\] and since by choice $N-M^{2}\leq 2M \leq 2N^{1/2}$, we have \[\left(\frac{N-M^{2}}{N}\right)^{\frac{s+1}{2^s}} \lesssim_{s} \left(\frac{1}{N^{1/2}}\right)^{\frac{s+1}{2^{s}}} \leq N^{-1/2^{s+2}},\] where we use that $(s+1)/2^{s+1} \geq 1/2^{s+2}$.
\end{proof}

The last step is to show that $Y_n$ doesn't affect pointwise behavior, which follows immediately from \Cref{p:PositiveResult} and \Cref{l:Ynuniform}.

\begin{cor}
\label{c:Ynpointwise}
    Almost surely in $(Y_n)$, the following holds: For every measure preserving system $(X,\mu,T)$ and $f\in L^{\infty}(\mu)$, one has \[\frac{1}{N}\sum_{n=1}^{N}Y_nf(T^nx) \to 0,\] pointwise almost everywhere.
\end{cor}

All elements are now set in place for us to prove \Cref{p:balayage}.

\begin{proof}[Proof of \Cref{p:balayage}]
    Fix $s\geq 2$ so that $\norm{a}_{U^s[N]}\to 0$. The three statements follow from the three lemmas we just proved. First, $d(E) = 1-\eps$ almost surely by \Cref{l:densityan}. Secondly, by \eqref{e:usefulexpressionYn}, we may split \[\norm{\one_E - (1-\eps)}_{U^s[N]} \leq \norm{Y_n}_{U^s[N]} + \eps\norm{a}_{U^s[N]},\] and thus, since by hypothesis $\norm{a}_{U^s[N]}\to 0$ and by \Cref{l:Ynuniform}, almost surely, $\norm{Y_n}_{U^s[N]}\to 0$ as well, we have \[\norm{\one_E - (1-\eps)}_{U^s[N]} \to 0\] with probability $1$. Finally, since $a$ is not pointwise good, there exists a measure preserving system $(X,\mu,T)$ and a function $f\in L^{\infty}(\mu)$ such that the averages \[\frac{1}{N}\sum_{n=1}^{N}a(n)f(T^nx)\] fail to converge on a set of positive measure. Writing
    \begin{align}
        \frac{1}{\abs{E\cap [N]}}\sum_{n\in E\cap [N]}f(T^nx) &= \frac{N}{\abs{E \cap [N]}}\cdot \frac{1}{N}\sum_{n=1}^{N}Y_nf(T^nx)\\
        &\qquad+ (1-\eps)\cdot \frac{N}{\abs{E \cap [N]}}\cdot \frac{1}{N}\sum_{n=1}^{N}f(T^nx)\\
        &\qquad\qquad+\eps\cdot\frac{N}{\abs{E \cap [N]}}\cdot\frac{1}{N}\sum_{n=1}^{N}a(n)f(T^nx)
    \end{align}
    and using that $N/\abs{E\cap [N]} \to 1/(1-\eps)$ almost surely, the first term is seen to almost surely converge pointwise almost everywhere by \Cref{c:Ynpointwise}, and by  Birkhoff's Ergodic Theorem so does the second one. Meanwhile, \[\eps\cdot \frac{N}{\abs{E \cap [N]}}\cdot \frac{1}{N}\sum_{n=1}^{N}a(n)f(T^nx)\] fails to converge on a set of positive measure with probability $1$, which implies that so does \[\frac{1}{\abs{E\cap [N]}}\sum_{n\in E\cap [N]}f(T^nx)\] with the same probability. Intersecting all three full probability events gets the result. 
\end{proof}

\begin{remark}
    Since the Gowers norms of $Y_n$ decay polynomially, they always beat the Gowers norms of the weights constructed by \Cref{t:MainResult}. This means that a stronger quantitative formulation of \Cref{t:GowersUniformBad} is also true; for every $\psi\colon \R\to \R_{+}$ with $\psi(x)\to \infty$ and $\psi(x)/x \to 0$ as $x \to \infty$ we can find a set $E\subseteq \N$ with $d(E)=1-\eps$ that is Gowers uniform with \[\norm{\one_{E}-d(E)}_{U^s[N]} \lesssim_{s} \left(\frac{\psi(\log N)}{\log N}\right)^{1/2^s}\] for every $s\geq2$ and sufficiently large $N\geq 1$, and yet $E$ is not pointwise good.
\end{remark}

\begin{remark}
    The same arguments as in \Cref{sec:Genericity} show that if for every $0 < \delta < 1$ we consider the space of subsets $E\subseteq \N$ with $d(E) > \delta$ and endow it with the metric given by the upper density of the symmetric difference \[d_{\triangle}(E,F) = \limsup_{N\to \infty}\frac{\abs{(E\triangle F) \cap [N]}}{N},\] then, after quotienting by $d_{\Delta}(E,F)=0$, the collection of subsets that are not pointwise good is again open and dense in the resulting space. In this scenario, it is crucial to restrict to subsets with densities uniformly bounded away from zero, since sets of $0$ density can have radically different averaging behavior, and in fact the natural analogue to \eqref{e:AveragesBesicovitch} does not hold. After this restriction, we do however have the inequality
    \begin{equation}
    \label{e:SubsequentialBesicovitch}
        \limsup_{N\to \infty} \norm{\frac{1}{\abs{E \cap [N]}}\sum_{n\in E\cap [N]}f(T^nx) - \frac{1}{\abs{F \cap [N]}}\sum_{n\in F\cap [N]}f(T^nx)} \leq \frac{2\norm{f}_{L^{\infty}(\mu)}}{\delta} d_{\triangle}(E,F),
    \end{equation}
    which allows the analogue of \Cref{l:Openness} to go through. The only step that needs care is the counterpart to \Cref{l:Suspension}. Here, instead of the suspension construction, given a pointwise good set $E$ and a Gowers uniform set $F$ that is not pointwise good, one can construct a random set $G$ using the probabilistic balayage technique with \[\P(n\in G) = (1-\eps)\one_{E} + \eps \one_{F}.\] Arguing as above shows that such a set $G$ will almost surely fail to be pointwise good, have density greater than $\delta$ and $d_{\Delta}(E,G)\leq 2\eps$. Furthermore, it will also preserve Gowers uniformity, so the genericity result will also hold even when restricting to Gowers uniform sets.
\end{remark}
\bigskip

This concludes the main body of the paper.
\appendix
\crefalias{section}{appendix}
\section{The Strong Sweeping-Out Property}
\label{appendix:Sweeping}

Most of the proofs and ideas laid out here are identical or simplifications of proofs given in the main body, so we only include details when strictly necessary. We begin by providing the proper context.\\

A measure preserving system $(X,\mu,T)$ on a standard probability space is said to be \textit{aperiodic} if the orbit $\set{T^nx \colon n\geq 0}$ of almost every point $x\in X$ is infinite. A weight $w\colon \Z\to \R$ is said to have the \textit{strong sweeping-out property} if for every aperiodic measure preserving system $(X,\mu,T)$ and every $\eps>0$, there exists a measurable set $E\subseteq X$ with $\mu(E)<\eps$, so that 
\[\liminf_{N\to\infty}\frac{1}{N}\sum_{n=1}^{N}w(n)\one_{E}(T^nx) = 0,\qquad \limsup_{N\to \infty}\frac{1}{N}\sum_{n=1}^{N}w(n)\one_{E}(T^nx)=1\] 
pointwise almost everywhere. This phenomenon was first encountered by U. Krengel \cite{Krengel} and then systematically studied by A. del Junco and J. Rosenblatt \cite{DelJuncoRosenblatt} and M. Akcoglu, Bellow, R.L. Jones, Losert, K. Reinhold-Larsson, and Wierdl \cite{AkcogluSweeping}. In particular, in \cite{AkcogluSweeping} they establish a general criterion for an averaging scheme to be strongly sweeping-out. We rephrase it here in the language of weighted ergodic averages; see also \cite[Proposition~2.11]{CespedesDonoso}.

\begin{theorem}[{\cite[Theorem~2.3]{AkcogluSweeping}}]
\label{t:SweepingCriterion}
    Let $w\colon \Z\to \R_{+}$ be a weight satisfying \[\frac{1}{N}\sum_{n=1}^{N}w(n) \to 1\] as $N\to \infty$. Suppose that for every $0<\eps<1$, every $M>0$ and every $N_0\geq 1$ there exists some $N_1\geq N_0$ and a finite set $E\subseteq \Z$ so that 
    \[\abs{\set{m\in \Z \colon \sup_{N_0\leq N \leq N_1} \frac{1}{N}\sum_{n=1}^{N}w(n)\one_{E}(m+n) > 1-\eps}} > M\cdot \abs{E}.\]
    Then $w$ is strongly sweeping-out.
\end{theorem}

Observe that a bounded weight will never be strongly sweeping-out as for ergodic systems one can always bound 
\[\limsup_{N\to \infty}\abs{\frac{1}{N}\sum_{n=1}^{N}w(n)\one_{E}(T^nx)} \leq \norm{w}_{\ell^{\infty}(\Z)}\cdot \mu(E)\] 
by Birkhoff's Ergodic Theorem. On the other hand, even for unbounded weights we have the following analogue of \Cref{t:PositiveResultFourier} and \Cref{p:PositiveResult}.

\begin{theorem}
\label{t:PositiveResultSweeping}
    Let $w\colon \Z\to \R_{+}$ be a weight admitting a density 
    \[d(w) = \lim_{N\to \infty}\frac{1}{N}\sum_{n=1}^{N}w(n)\]
    and suppose that at least one of the following conditions holds:
    \begin{enumerate}
        \item There exists some $\eta>1/2$ such that for every $N\geq 2$
        \[\sup_{\beta\in \T}\abs{\frac{1}{N}\sum_{n=1}^{N}(w(n)-d(w))e(n\beta)}\lesssim \left(\frac{1}{\log N}\right)^{\eta}.\]
        \item There exists some $s\geq 2$ and some $\eta > 1/2^s$ such that for every $N\geq 2$
        \[\norm{w-d(w)}_{U^s[N]}\lesssim_{s} \left(\frac{1}{\log N}\right)^{\eta}.\]
    \end{enumerate}
    Then $w$ is pointwise good for $L^{\infty}$ functions.
\end{theorem}

\begin{remark}
    We do not need to assume that $d(w)$ exists a priori; if either hypothesis is true for $w-c$ where $c\in \R$ is some constant, then $w$ immediately admits a density and $d(w)=c$.
\end{remark}

\begin{proof}
    Let $(X,\mu,T)$ be a measure preserving system, $f\in L^{\infty}(\mu)$ and fix $\rho>1$. By linearity, we may assume that $f\geq 0$. Proceeding as in \Cref{p:PositiveResult} shows that whenever one of the conditions above holds, then
    \[\frac{1}{\lfloor \rho^m \rfloor}\sum_{n=1}^{\lfloor \rho^m \rfloor}(w(n)-d(w))f(T^nx) \to 0,\]
    pointwise almost everywhere as $m\to \infty$. Birkhoff's Ergodic Theorem then implies that
    \[\frac{1}{\lfloor \rho^m \rfloor}\sum_{n=1}^{\lfloor \rho^m \rfloor}w(n)f(T^nx) \to d(w)P_{\II}f(x),\]
    where $P_{\II}f$ denotes the projection of $f$ onto the subspace of $L^2(\mu)$ of $T$-invariant functions. Since $w$ and $f$ are non-negative, its partial sums are non-decreasing, so for $\lfloor \rho^m \rfloor \leq N < \lfloor \rho^{m+1} \rfloor$ we may pointwise bound
    \[\frac{\lfloor \rho^m \rfloor}{\lfloor \rho^{m+1} \rfloor}\cdot \frac{1}{\lfloor \rho^m \rfloor}\sum_{n=1}^{\lfloor \rho^m \rfloor}w(n)f(T^nx) \leq \frac{1}{N}\sum_{n=1}^{N}w(n)f(T^nx) \leq \frac{\lfloor \rho^{m+1} \rfloor}{\lfloor \rho^m \rfloor}\cdot \frac{1}{\lfloor \rho^{m+1} \rfloor}\sum_{n=1}^{\lfloor \rho^{m+1} \rfloor}w(n)f(T^nx),\]
    and letting $N\to \infty$ and hence $m\to \infty$ implies
    \[\frac{1}{\rho}d(w)P_{\II}f(x) \leq \liminf_{N\to \infty}\frac{1}{N}\sum_{n=1}^{N}w(n)f(T^nx) \leq \limsup_{N\to \infty} \frac{1}{N}\sum_{n=1}^{N}w(n)f(T^nx) \leq \rho d(w)P_{\II}f(x).\]
    Letting $\rho\searrow 1$ in a countable fashion implies that
    \[\frac{1}{N}\sum_{n=1}^{N}w(n)f(T^nx) \to d(w)P_{\II}f(x)\] pointwise almost everywhere.
\end{proof}

In particular, any such weight can never be strong sweeping-out. The divergence mechanism underlying the construction in \Cref{sec:Construction} can be modified to yield the following analogue to \Cref{t:MainResult}.

\begin{theorem}
\label{t:SweepingOut}
    For any pair of non-decreasing functions $\psi,\Psi\colon \R_{+}\to \R_{+}$ with $\psi(x),\Psi(x) \to \infty$ as $x\to \infty$, $\Psi\geq 1$, there exists a strongly sweeping-out weight $u\colon \Z\to \Z_{+}$ satisfying
    \[\lim_{N\to \infty} \frac{1}{N}\sum_{n=1}^{N}u(n)=1, \qquad u(n)\leq \Psi(n),\ \forall n\geq 1,\]
    so that the centered weight $v=u-1$ satisfies, for every $N\geq 2$,
    \[\sup_{\beta \in \T}\abs{\frac{1}{N}\sum_{n=1}^{N}v(n)e(n\beta)} \lesssim \left(\frac{\psi(\log N)}{\log N}\right)^{1/2}\]
    and 
    \[\norm{v}_{U^s[N]}\lesssim_{s} \left(\frac{\psi(\log N)}{\log N}\right)^{1/2^s},\]
    simultaneously for every $s\geq 2$.
\end{theorem}

In particular, if $\psi$ is chosen so that the quotient above decays, we obtain the following consequence in line with \Cref{c:StronglyErgodicNotPointwiseGood}.

\begin{cor}
    There exists a non-negative, integer valued and strongly ergodic weight that satisfies the strong sweeping-out property.
\end{cor}

Before outlining the construction, recall from \Cref{subsec:Multiplicative}, that given $q\geq 2$ and a prime $p$ such that $q\mid p-1$ there always exists a multiplicative character of $\F_{p}$ of order $q$. Denote said character by $\chi_{p,q}$ and observe that by definition 
\begin{equation}
\label{e:supportsize}
    \abs{\set{t\in \F_{p} \colon \chi_{p,q}(t)=1}} = \frac{p-1}{q}.
\end{equation}
Starting from $\chi_{p,q}$, we define a function $u_{p,q}\colon \F_{p} \to \set{0,1,q}$ by 
\[u_{p,q}(t) = 1 + \sum_{r=1}^{q-1}\chi_{p,q}(t)^{r} = 
\begin{cases}
    1, &\text{ if } t=0,\\
    0, &\text{ if } t\neq 0, \chi_{p,q}(t)\neq 1,\\
    q, &\text{ if } t\neq 0, \chi_{p,q}(t)= 1,\\
\end{cases}\]
and let $v_{p,q} = u_{p,q}-1$. Observe that by \eqref{e:cancellationidentity}, we have
\begin{equation}
\label{e:sumupqvpq}
    \sum_{t\in \F_{p}}v_{p,q}(t)=0, \qquad \sum_{t\in \F_{p}}u_{p,q}(t)=p.
\end{equation}

Furthermore, we have the following estimates.
 
\begin{lemma}
    For any $q\geq 2$ and prime $p$ with $q\mid p-1$ we have
    \[\max_{a\in \F_{p}} \abs{\widehat{v_{p,q}}(a)} \leq \frac{q}{p^{1/2}}\]
    and 
    \[\norm{v_{p,q}}_{U^s(\Z_{p})} \lesssim_{s} \frac{q}{p^{1/2^s}}\]
    for every $s\geq 2$.
\end{lemma}

The Fourier estimate follows directly from the triangle inequality and the Gauss-sum estimate given in \Cref{l:GaussSum}, after noticing that the maps $t\mapsto \chi_{p,q}(t)^{r}$ are also multiplicative characters of $\F_p$ for $1\leq r < q$. The second statement is proven in the same way and only requires the following consequence of Weil's Bound that says that cyclic Gowers norms of non-trivial multiplicative characters exhibit cancellation at the natural exponent. The proof is essentially a simplified version of that of \Cref{l:CharacterWordGowers}.

\begin{lemma}
\label{l:GowersMultiplicative}
    Let $\chi$ be a non-trivial multiplicative character of $\F_{p}$ and let $s\geq 2$. Then 
    \[\norm{\chi}_{U^s(\Z_p)}^{2^s} \lesssim_{s} \frac{1}{p}.\]
\end{lemma}

Extending each $v_{p,q}$ $p$-periodically to all of $\Z$ and arguing as in \Cref{l:FourierDecay} and \Cref{c:GowersNormsCharacterWord} one can transfer these estimates to $\Z$ and obtain the following.

\begin{lemma}
\label{l:vpqcancellation}
    Let $q\geq 2$, $p\geq 3$ be a prime with $q\mid p-1$. Then, for any $N\geq 1$, any interval $I\subseteq [N]$ and any shift $a\in \Z$ one has 
    \[\sup_{\beta\in \T}\abs{\sum_{n\in I}v_{p,q}(n+a)e(n\beta)} \lesssim \frac{q}{p^{1/2}}\left(\abs{I}+p \log p\right)\]
    and if $\abs{I}\geq p^{2}$, then for every $s\geq 2$ 
    \[\norm{v_{p,q}(\cdot +a)}_{U^s[I]}\lesssim \frac{q}{p^{1/2^s}}.\]
\end{lemma}

The construction is now very similar to what we did in \Cref{sec:Construction}. Fix $\psi$ and $\Psi$ as in \Cref{t:SweepingOut}, put $q_k = k+2$ for every $k\geq 0$ and set $B_0=1$. At stage $k\geq 1$, by Dirichlet's Theorem on arithmetic progressions, we can find a prime $p_k$ with $p_k \equiv 1 \mod q_k$ sufficiently large so that the following conditions hold simultaneously 
\[p_k > B_{k-1},\quad \psi(\log p_k) \geq q_k^{2^k}\log q_k, \quad \Psi(p_k^{4}) \geq q_k,\]
and define 
\[A_k \coloneqq p_k^{4},\quad B_{k} \coloneqq p_{k}^{4}q_k^{p_k},\]
along with the interval 
\[I^{(k)} \coloneqq (A_k, B_k] \cap \Z,\]
subdivided into the disjoint blocks 
\[I^{(k)}_{j} \coloneqq (A_kq_k^{j-1},A_{k}q_{k}^{j}] \cap \Z, \quad \abs{I^{(k)}_{j}} = A_kq_{k}^{j-1}(q_k-1)\]
for $1\leq j \leq p_k$. Extending all $u_{p_k,q_k}$ $p_k$-periodically to all of $\Z$, we can then define a weight $u\colon \Z\to \Z$ via
\begin{equation}
    u(n) = 
    \begin{cases}
        u_{p_k,q_k}(n+j), &\text{ if } n\in I_j^{(k)}, \text{ for some }k\geq 1 \text{ and }1\leq j\leq p_k,\\
        1, &\text{ otherwise,}
    \end{cases}
\end{equation}
and put 
\begin{equation}
    v(n)= u(n)-1 = 
    \begin{cases}
        v_{p_k,q_k}(n+j), &\text{ if } n\in I_j^{(k)}, \text{ for some }k\geq 1 \text{ and }1\leq j\leq p_k,\\
        0, &\text{ otherwise.}
    \end{cases}
\end{equation}

Observe that in the $k$-th stage, $u$ takes the values $0$, $1$ and $q_k$. It follows that $u\geq 0$ and that $u(n)\leq \Psi(n)$. Moreover, by \eqref{e:sumupqvpq} and the fact that the length of each $I_{j}^{(k)}$ is a multiple of $p_k$, we have 
\begin{equation}
\label{e:sumuv}
    \sum_{n\in I_{j}^{(k)}}v(n) = 0, \qquad \sum_{n\in I^{(k)}_{j}}u(n) = \abs{I_{j}^{(k)}},
\end{equation}
for every $k\geq 1$ and $1\leq j \leq p_k$. Since $v\equiv 0$ outside of the blocks $I^{(k)}$, 
\[\abs{\sum_{n=1}^{N}u(n) - N} = \abs{\sum_{n=1}^{N}v(n)} \leq q_kp_k,\]
whenever $A_k\leq N\leq B_k$ and hence 
\[\frac{1}{N}\sum_{n=1}^{N}u(n) \to 1\]
as $N\to \infty$. Using \Cref{l:vpqcancellation} and arguing exactly as in \Cref{l:FourierCancellation} and \Cref{l:GowersCancellation}, we obtain the bounds 
\[\sup_{\beta\in \T}\abs{\frac{1}{N}\sum_{n=1}^{N}v(n)e(n\beta)} \lesssim \frac{q_k}{p_k^{1/2}},\quad \norm{v}_{U^s[N]}\lesssim_{s} \frac{q_k}{p_k^{1/2^{s}}},\ \forall s\geq 2,\]
whenever $A_k\leq N \leq B_k$; appealing to our choice of parameters as in \Cref{p:ConcreteRate} yields 
\[\sup_{\beta \in \T}\abs{\frac{1}{N}\sum_{n=1}^{N}v(n)e(n\beta)} \lesssim \left(\frac{\psi(\log N)}{\log N}\right)^{1/2}\]
and 
\[\norm{v}_{U^s[N]}\lesssim_{s} \left(\frac{\psi(\log N)}{\log N}\right)^{1/2^s}\]
simultaneously for every $s\geq 2$. So to establish \Cref{t:SweepingOut}, it only remains to show that $u$ is strongly sweeping-out. 

\begin{proof}[Proof of \Cref{t:SweepingOut}]
    We use \Cref{t:SweepingCriterion} and begin by laying out what happens at each stage $k\geq 1$. Regarding $u_{p_k,q_k}$ as a function on $\F_{p_k}$ and letting 
    \[E_k \coloneqq \set{m \in \Z \colon m \in \supp u_{p_k,q_k}  \mod p_{k}} \cap \set{1,\ldots, 2B_k},\]
    one can observe that by \eqref{e:supportsize} and the fact that $p_k \mid B_k$, we have 
    \[\abs{E_k} \leq \frac{2B_k}{p_k}\left(\frac{p_k-1}{q_k} + 1\right) < \frac{4B_k}{q_k}.\]
    Now given $1\leq m\leq B_k$, if we let $1\leq j\leq p_k$ such that $j \equiv m \mod p_k$, then, for $n\in I_{j}^{(k)}$, since $1\leq n+m \leq 2B_k$, we have 
    \[u(n)=u_{p_k,q_k}(n+j) = u(n)\one_{E_k}(n+m).\]
    In particular, since $u\geq 0$, by \eqref{e:sumuv} and our choice of parameters 
    \[\frac{1}{A_kq_k^{j}}\sum_{n=1}^{A_kq_k^{j}}u(n)\one_{E_k}(n+m) \geq \frac{1}{A_kq_k^{j}}\sum_{n\in I_{j}^{(k)}}u(n) = \frac{1}{A_kq_k^{j}}\abs{I_j^{(k)}} = 1 - \frac{1}{q_k},\]
    so for every $k\geq 1$ and every $1\leq m \leq B_k$ 
    \[\sup_{A_k\leq N\leq B_k} \frac{1}{N}\sum_{n=1}^{N}u(n)\one_{E_k}(n+m) \geq 1- \frac{1}{q_k}.\]
    This means that given any $0<\eps<1$, any $M>0$, and any $N_0\geq 1$, by choosing $k\geq 1$ large enough so that \[A_k \geq N_0, \quad \frac{1}{q_k} < \eps, \quad q_k > 4M,\] and taking $N_1=B_k$, we have 
    \[\abs{\set{m\in \Z \colon \sup_{N_0\leq N \leq N_1} \frac{1}{N}\sum_{n=1}^{N}u(n)\one_{E_k}(m+n) > 1-\eps}} \geq B_k,\]
    since every $1\leq m\leq B_k$ will belong to this set. On the other hand, we also have 
    \[M\abs{E_k} < \frac{4MB_k}{q_k} < B_k,\]
    so putting both things together we see that the criterion of \Cref{t:SweepingCriterion} is satisfied and $u$ is strongly sweeping-out.
\end{proof}

\printbibliography

@article{AkcogluSweeping,
  title={{The Strong Sweeping Out Property for Lacunary Sequences, {Riemann} Sums, Convolution Powers, and Related Matters}},
  author={Akcoglu, Mustafa and Bellow, Alexandra and Jones, Roger L. and Losert, Viktor and Reinhold-Larsson, Karin and Wierdl, M{\'a}t{\'e}},
  journal={Ergodic Theory and Dynamical Systems},
  volume={16},
  number={2},
  pages={207--253},
  year={1996},
  publisher={Cambridge University Press},
  doi={10.1017/S0143385700008798}
}

@article{BellowLosert,
  title={{The Weighted Pointwise Ergodic Theorem and the Individual Ergodic Theorem Along Subsequences}},
  author={Bellow, Alexandra and Losert, Viktor},
  journal={Transactions of the American Mathematical Society},
  volume={288},
  number={1},
  pages={307--345},
  year={1985},
  publisher={American Mathematical Society},
  doi={10.1090/S0002-9947-1985-0773063-8}
}

@article{Birkhoff,
  title={{Proof of the Ergodic Theorem}},
  author={Birkhoff, George D.},
  journal={Proceedings of the National Academy of Sciences of the United States of America},
  volume={17},
  number={12},
  pages={656--660},
  year={1931},
  publisher={National Academy of Sciences},
  doi={10.1073/pnas.17.2.656}
}

@incollection{BourgainPrimes,
  title={{An Approach to Pointwise Ergodic Theorems}},
  author={Bourgain, Jean},
  booktitle={{Geometric Aspects of Functional Analysis: Israel Seminar (GAFA) 1986--87}},
  editor={Lindenstrauss, J. and Milman, V. D.},
  series={Lecture Notes in Mathematics},
  volume={1317},
  pages={204--223},
  year={1988},
  publisher={Springer-Verlag},
  address={Berlin},
  doi={10.1007/BFb0081742}
}

@article{B0,
  title={{On the Maximal Ergodic Theorem for Certain Subsets of the Integers}},
  author={Bourgain, Jean},
  journal={Israel Journal of Mathematics},
  volume={61},
  number={1},
  pages={39--72},
  year={1988},
  doi={10.1007/BF02776301}
}

@article{BourgainReturnTimes,
  title={{Temps de Retour pour les Syst{\`e}mes Dynamiques}},
  author={Bourgain, Jean},
  journal={C. R. Acad. Sci. Paris S{\'e}r. I Math.},
  volume={306},
  number={12},
  pages={483--485},
  year={1988},
  mrnumber={939436}
}

@techreport{BourgainUnpublished,
  title={{Return Times of Dynamical Systems}},
  author={Bourgain, Jean},
  type={Preprint},
  number={IHES/M/88/15},
  year={1988},
  institution={Institut des Hautes {\'E}tudes Scientifiques},
  url={https://repo-archives.ihes.fr/FONDS_IHES/I_Prepublications/BOURGAIN/1984-1989/M_88_15/M_88_15_web.pdf}
}

@article{BFKO,
  title={{Appendix on Return Times of Dynamical Systems}},
  author={Bourgain, Jean and Furstenberg, Harry and Katznelson, Yitzhak and Ornstein, Donald S.},
  journal={Publications Math{\'e}matiques de l'Institut des Hautes {\'E}tudes Scientifiques},
  volume={69},
  pages={42--45},
  year={1989},
  publisher={Institut des Hautes {\'E}tudes Scientifiques},
  note={Appendix to Jean Bourgain's \emph{Pointwise Ergodic Theorems for Arithmetic Sets}},
  doi={10.1007/BF02698839}
}

@article{C1,
  title={{Ergodic Theory and Translation-Invariant Operators}},
  author={Calder{\'o}n, Alberto P.},
  journal={Proceedings of the National Academy of Sciences of the United States of America},
  volume={59},
  number={2},
  pages={349--353},
  year={1968},
  publisher={National Academy of Sciences},
  doi={10.1073/pnas.59.2.349}
}

@article{CespedesDonoso,
  title={{Pointwise Ergodic Averages Along the Omega Function in Number Fields}},
  author={C{\'e}spedes, Diego and Donoso, Sebasti{\'a}n},
  journal={Bulletin of the London Mathematical Society},
  volume={58},
  number={8},
  pages={e70450},
  year={2026},
  publisher={John Wiley \& Sons},
  doi={10.1112/blms.70450}
}

@article{Christ,
  title={{A Weak Type (1,1) Inequality for Maximal Averages Over Certain Sparse Sequences}},
  author={Christ, Michael},
  journal={arXiv preprint arXiv:1108.5664},
  year={2011},
  doi={10.48550/arXiv.1108.5664}
}

@article{DelJuncoRosenblatt,
  title={{Counterexamples in Ergodic Theory and Number Theory}},
  author={del Junco, Andr{\'e}s and Rosenblatt, Joseph},
  journal={Mathematische Annalen},
  volume={245},
  number={3},
  pages={185--197},
  year={1979},
  publisher={Springer},
  doi={10.1007/BF01673506}
}

@book{EisnerFarkas,
  title={{A Journey Through Ergodic Theorems}},
  author={Eisner, Tanja and Farkas, B{\'a}lint},
  series={Birkh{\"a}user Advanced Texts Basler Lehrb{\"u}cher},
  year={2025},
  publisher={Birkh{\"a}user},
  address={Cham},
  doi={10.1007/978-3-031-96506-7}
}

@article{EKLR,
  title={{The {Chowla} and the {Sarnak} Conjectures From Ergodic Theory Point of View}},
  author={El Abdalaoui, El Houcein and Ku{\l}aga-Przymus, Joanna and Lema{\'n}czyk, Mariusz and de la Rue, Thierry},
  journal={Discrete and Continuous Dynamical Systems},
  volume={37},
  number={6},
  pages={2899--2944},
  year={2017},
  doi={10.3934/dcds.2017125}
}

@article{Fan,
  title={{Weighted {Birkhoff} Ergodic Theorem With Oscillating Weights}},
  author={Fan, Ai-Hua},
  journal={Ergodic Theory and Dynamical Systems},
  volume={39},
  number={5},
  pages={1275--1289},
  year={2019},
  publisher={Cambridge University Press},
  doi={10.1017/etds.2017.81}
}

@article{FanKonieczny,
  title={{On Uniformity of {$q$}-Multiplicative Sequences}},
  author={Fan, Ai-Hua and Konieczny, Jakub},
  journal={Bulletin of the London Mathematical Society},
  volume={51},
  number={3},
  pages={466--488},
  year={2019},
  publisher={John Wiley \& Sons},
  doi={10.1112/blms.12245}
}

@article{FHK2007,
  title={{Multiple Recurrence and Convergence for Sequences Related to the Prime Numbers}},
  author={Frantzikinakis, Nikos and Host, Bernard and Kra, Bryna},
  journal={Journal f{\"u}r die reine und angewandte Mathematik},
  volume={2007},
  number={611},
  pages={131--144},
  year={2007},
  doi={10.1515/CRELLE.2007.076}
}

@article{FHK2013,
  title={{The Polynomial Multidimensional {Szemer{\'e}di} Theorem Along Shifted Primes}},
  author={Frantzikinakis, Nikos and Host, Bernard and Kra, Bryna},
  journal={Israel Journal of Mathematics},
  volume={194},
  number={1},
  pages={331--348},
  year={2013},
  publisher={Springer},
  doi={10.1007/s11856-012-0132-y}
}

@article{FornalKrause2026,
  title={{Double Recurrence and Almost Sure Convergence: Primes and Weighted Theory}},
  author={Fornal, Jan and Krause, Ben},
  journal={arXiv preprint arXiv:2603.22203},
  year={2026},
  doi={10.48550/arXiv.2603.22203}
}

@article{FornalWW,
  title={{The {Wiener Wintner} Theorem Along the Primes}},
  author={Fornal, Jan and Fragkos, Anastasios and Krause, Ben and Lacey, Michael and Mousavi, Hamed and Sun, Yu-Chen},
  journal={arXiv preprint arXiv:2601.10459},
  year={2026},
  doi={10.48550/arXiv.2601.10459}
}

@article{FrantzikinakisHost2017,
  title={{Multiple Ergodic Theorems for Arithmetic Sets}},
  author={Frantzikinakis, Nikos and Host, Bernard},
  journal={Transactions of the American Mathematical Society},
  volume={369},
  number={10},
  pages={7085--7105},
  year={2017},
  doi={10.1090/tran/6870}
}

@article{GreenTao,
  title={{Linear Equations in Primes}},
  author={Green, Ben and Tao, Terence},
  journal={Annals of Mathematics},
  volume={171},
  number={3},
  pages={1753--1850},
  year={2010},
  publisher={Princeton University and Institute for Advanced Study},
  doi={10.4007/annals.2010.171.1753}
}

@article{HernandezRadic,
  title={{Infinite Prime Sumsets in Structured and {$U^{k}(\Phi)$}-Uniform Sets}},
  author={Hern{\'a}ndez, Felipe and Radi{\'c}, Trist{\'a}n},
  journal={arXiv preprint arXiv:2609.20522},
  year={2026},
  doi={10.48550/arXiv.2609.20522}
}

@article{HostKra,
  title={{Uniformity Seminorms on {$\ell^{\infty}$} and Applications}},
  author={Host, Bernard and Kra, Bryna},
  journal={Journal d'Analyse Math{\'e}matique},
  volume={108},
  pages={219--276},
  year={2009},
  doi={10.1007/s11854-009-0024-1}
}

@book{IwaniecKowalski,
  title={{Analytic Number Theory}},
  author={Iwaniec, Henryk and Kowalski, Emmanuel},
  series={Colloquium Publications},
  volume={53},
  year={2004},
  publisher={American Mathematical Society},
  address={Providence, RI},
  doi={10.1090/coll/053}
}

@article{Keane,
  title={{Generalized {Morse} Sequences}},
  author={Keane, Michael},
  journal={Zeitschrift f{\"u}r Wahrscheinlichkeitstheorie und Verwandte Gebiete},
  volume={10},
  number={4},
  pages={335--353},
  year={1968},
  publisher={Springer},
  doi={10.1007/BF00531855}
}

@article{Konieczny,
  title={{Gowers Norms for the {Thue--Morse} and {Rudin--Shapiro} Sequences}},
  author={Konieczny, Jakub},
  journal={Annales de l'Institut Fourier},
  volume={69},
  number={4},
  pages={1897--1913},
  year={2019},
  doi={10.5802/aif.3285}
}

@article{Koutsogiannis,
  title={{Closest Integer Polynomial Multiple Recurrence Along Shifted Primes}},
  author={Koutsogiannis, Andreas},
  journal={Ergodic Theory and Dynamical Systems},
  volume={38},
  number={2},
  pages={666--685},
  year={2018},
  publisher={Cambridge University Press},
  doi={10.1017/etds.2016.40}
}

@article{Krengel,
  title={{On the Individual Ergodic Theorem for Subsequences}},
  author={Krengel, Ulrich},
  journal={The Annals of Mathematical Statistics},
  volume={42},
  number={3},
  pages={1091--1095},
  year={1971},
  publisher={Institute of Mathematical Statistics},
  doi={10.1214/aoms/1177693338}
}

@article{Leibman,
  title={{Multiple Recurrence Theorem for Measure Preserving Actions of a Nilpotent Group}},
  author={Leibman, A.},
  journal={Geometric and Functional Analysis},
  volume={8},
  number={5},
  pages={853--931},
  year={1998},
  publisher={Springer},
  doi={10.1007/s000390050077}
}

@article{LesignePolynomial,
  title={{Spectre Quasi-Discret et Th{\'e}or{\`e}me Ergodique de {Wiener-Wintner} pour les Polyn{\^o}mes}},
  author={Lesigne, Emmanuel},
  journal={Ergodic Theory Dynam. Systems},
  volume={13},
  number={4},
  pages={767--784},
  year={1993},
  mrnumber={1257033}
}

@article{LesigneQuasRosenblattWierdl,
  title={{Generation of Measures on the Torus With Good Sequences of Integers}},
  author={Lesigne, Emmanuel and Quas, Anthony and Rosenblatt, Joseph and Wierdl, M{\'a}t{\'e}},
  journal={Indiana University Mathematics Journal},
  volume={74},
  number={6},
  pages={1801--1847},
  year={2025},
  doi={10.1512/IUMJ.2025.74.60561}
}

@article{LinOlsenTempelman,
  title={{On Modulated Ergodic Theorems for {Dunford-Schwartz} Operators}},
  author={Lin, Michael and Olsen, James and Tempelman, Arkady},
  journal={Illinois Journal of Mathematics},
  volume={43},
  number={3},
  pages={542--567},
  year={1999},
  publisher={University of Illinois},
  doi={10.1215/ijm/1255985110}
}

@article{MatthiesenTeravainenWang,
  title={{Quantitative Asymptotics for Polynomial Patterns in the Primes}},
  author={Matthiesen, Lilian and Ter{\"a}v{\"a}inen, Joni and Wang, Mengdi},
  journal={Mathematika},
  volume={72},
  number={3},
  pages={e70103},
  year={2026},
  publisher={London Mathematical Society},
  doi={10.1112/mtk.70103}
}

@article{Peluse,
  title={{Bounds for Sets With No Polynomial Progressions}},
  author={Peluse, Sarah},
  journal={Forum of Mathematics, Pi},
  volume={8},
  pages={e16},
  year={2020},
  publisher={Cambridge University Press},
  doi={10.1017/fmp.2020.11}
}

@article{Radic,
  title={{Infinite Sumsets in {$U^{k}(\Phi)$}-Uniform Sets}},
  author={Radi{\'c}, Trist{\'a}n},
  journal={arXiv preprint arXiv:2601.06915},
  year={2026},
  doi={10.48550/arXiv.2601.06915}
}

@article{Rudolph,
  title={{Fully Generic Sequences and a Multiple-Term Return-Times Theorem}},
  author={Rudolph, Daniel J.},
  journal={Inventiones mathematicae},
  volume={131},
  number={1},
  pages={199--228},
  year={1998},
  publisher={Springer},
  doi={10.1007/s002220050202}
}

@incollection{RyllNardzewski,
  title={{Topics in Ergodic Theory}},
  author={Ryll-Nardzewski, Czes{\l}aw},
  booktitle={{Probability-Winter School: Proceedings of the Fourth Winter School on Probability Held at Karpacz, Poland, January 1975}},
  editor={Ciesielski, Z. and Urbanik, K. and Woyczy{\'n}ski, W. A.},
  series={Lecture Notes in Mathematics},
  volume={472},
  pages={131--156},
  year={1975},
  publisher={Springer-Verlag},
  address={Berlin},
  doi={10.1007/BFb0081951}
}

@article{Tal,
  title={{The Return Times Theorem, Auto-Correlation and Sequences With an Empty {Fourier--Bohr} Spectrum}},
  author={Tal, Matan},
  journal={Studia Mathematica},
  volume={288},
  number={2},
  pages={135--154},
  year={2026},
  doi={10.4064/sm250117-8-1}
}

@article{TaoZiegler,
  title={{Polynomial Patterns in the Primes}},
  author={Tao, Terence and Ziegler, Tamar},
  journal={Forum of Mathematics, Pi},
  volume={6},
  pages={e1},
  year={2018},
  publisher={Cambridge University Press},
  doi={10.1017/fmp.2017.3}
}

@article{Teravainen,
  title={{Pointwise Convergence of Ergodic Averages With {M{\"o}bius} Weight}},
  author={Ter{\"a}v{\"a}inen, Joni},
  journal={arXiv preprint arXiv:2401.03174},
  year={2024},
  doi={10.48550/arXiv.2401.03174}
}

@unpublished{ThouvenotWeiss,
  title={{Mixing Sequences That Are Not Self Sampling}},
  author={Thouvenot, Jean-Paul and Weiss, Benjamin},
  note={Privately communicated notes}
}

@article{VonNeumman,
  title={{Proof of the Quasi-Ergodic Hypothesis}},
  author={von Neumann, John},
  journal={Proceedings of the National Academy of Sciences of the United States of America},
  volume={18},
  number={1},
  pages={70--82},
  year={1932},
  publisher={National Academy of Sciences},
  doi={10.1073/pnas.18.1.70}
}

@article{Walsh,
  title={{Norm Convergence of Nilpotent Ergodic Averages}},
  author={Walsh, Miguel N.},
  journal={Annals of Mathematics},
  volume={175},
  number={3},
  pages={1667--1688},
  year={2012},
  publisher={Princeton University and Institute for Advanced Study},
  doi={10.4007/annals.2012.175.3.15}
}

@article{WienerWintner,
  title={{Harmonic Analysis and Ergodic Theory}},
  author={Wiener, Norbert and Wintner, Aurel},
  journal={American Journal of Mathematics},
  volume={63},
  number={2},
  pages={415--426},
  year={1941},
  publisher={Johns Hopkins University Press},
  doi={10.2307/2371534}
}

@article{WierdlPrimes,
  title={{Pointwise Ergodic Theorem Along the Prime Numbers}},
  author={Wierdl, M{\'a}t{\'e}},
  journal={Israel Journal of Mathematics},
  volume={64},
  number={3},
  pages={315--336},
  year={1988},
  publisher={Springer},
  doi={10.1007/BF02882425}
}

@article{ZorinKranichDouble,
  title={{A Double Return Times Theorem}},
  author={Zorin-Kranich, Pavel},
  journal={Israel Journal of Mathematics},
  volume={229},
  number={1},
  pages={255--267},
  year={2019},
  publisher={Springer},
  doi={10.1007/s11856-018-1797-7}
}

\end{document}